\documentclass[11pt]{amsart}
\usepackage[margin=1.2in]{geometry}
\usepackage[utf8]{inputenc}
\usepackage{graphicx}
\usepackage{amssymb, amsmath, amsthm, graphics}
\usepackage{mathabx,epsfig}
\usepackage[mathscr]{euscript}
\usepackage{tikz}
\usetikzlibrary{calc}
\usetikzlibrary{matrix}
\usepackage{appendix}

\usepackage{latexsym}
\usepackage[all]{xy}
\usepackage{color}

\usepackage{bbold}

\usepackage{tikz}
\usepackage{tikz-cd}
\input xy
\xyoption{all}
\usepackage{euscript}
\usepackage{multirow}
\usepackage{etex, pictexwd,dcpic}
\usetikzlibrary{positioning}
\usetikzlibrary{shapes.geometric}
\usetikzlibrary{shapes.misc}
\usetikzlibrary{calc}
\usetikzlibrary{positioning}

\usepackage{pgflibraryarrows}
\usepackage{pgflibrarysnakes}

\usetikzlibrary{trees} 
\usetikzlibrary[trees] 

\theoremstyle{plain}
\newtheorem{theorem}{Theorem}[section]

\newtheorem{proposition}[theorem]{Proposition}

\newtheorem{lemma}[theorem]{Lemma}
\theoremstyle{definition}
\newtheorem{definition}[theorem]{Definition}
\newtheorem{notation}[theorem]{Notation}

\theoremstyle{remark}

\numberwithin{equation}{section}

\renewcommand{\(}{\begin{equation*}}
\renewcommand{\)}{\end{equation*}}
\newcommand{\bea}{\begin{eqnarray*}}
\newcommand{\eea}{\end{eqnarray*}}

\def\proof {{Proof.}\hspace{7pt}}

\def\endofproof {\hfill{$\Box$}\\}

\newcommand{\beq}{\begin{equation}}
\newcommand{\eeq}{\end{equation}}

\numberwithin{equation}{section}

\renewcommand{\(}{\begin{equation}}
\renewcommand{\)}{\end{equation}}

\def\1{{\bf 1}}

\def\<{\langle}
\def\>{\rangle}

\numberwithin{equation}{section}


\newcommand{\NN}{\ensuremath{\mathbb N}}

\newcommand{\ZZ}{\ensuremath{\mathbb Z}}
\newcommand{\Z}{\ensuremath{\mathbb Z}} 

\newcommand{\HH}{\ensuremath{\mathbb H}}

\makeatletter
\newcommand{\colim@}[2]{%
  \vtop{\m@th\ialign{##\cr
    \hfil$#1\operator@font colim$\hfil\cr
    \noalign{\nointerlineskip\kern1.5\ex@}#2\cr
    \noalign{\nointerlineskip\kern-\ex@}\cr}}%
}
\newcommand{\colim}{%
  \mathop{\mathpalette\colim@{\rightarrowfill@\textstyle}}\nmlimits@
}

\makeatother

\DeclareRobustCommand{\SkipTocEntry}[5]{}

\title{Cohomotopy and flux quantization in $M$-theory}
\date{\today}
\address{Department of Mathematics, Statistics, and Physics, Wichita State University, 1845 N. Fairmount, Wichita, KS 67260}
\author{Daniel Grady}

\begin{document}

\begin{abstract}
We identify some of the $k$-invariants for the Postnikov tower of the stable and unstable 4-sphere. Assuming the stable Hypothesis H of Fiorenza--Sati--Schreiber, we use the resulting obstruction theory to prove that the Chern--Simons term in the effective action of M-theory is well defined. In particular, we do not assume the presence of an $E_8$-gauge field.
\end{abstract}

\maketitle

\tableofcontents

\section{Introduction}

In this note, we follow up on several remarks made in \cite{GS}. We give a self-contained identification of some k-invariants of the Postnikov tower for both the stable and unstable 4-sphere. We use the obstruction theory to prove that for a closed oriented 12-manifold $Q$ and a de Rham class $G\in H^4_{dR}(M)$ with integral periods, the cubic term
\begin{equation}\label{witint} \frac{1}{6}\int_{Q}G\wedge G \wedge G 
\end{equation}
is integral, provided that $G$ lifts to a stable cohomotopy class in degree $4$. 

The integrality of the above expression is necessary to get a well defined action functional for $11$D supergravity, in its usual normalization (see e.g., \cite[page 6]{Fig}). Traditionally, this problem has received attention after including an expected quantum correction (a ``1-loop term") and arguments for the integrality of \eqref{witint} using $E_8$-gauge theory were given in \cite{W}. However, since there is no $E_8$-gauge field present in supergravity, one expects to obtain integrality without reference to such a gauge field. This problem was recognized in \cite{W}, where it was suggested that the integrality of \eqref{witint} should hold without assuming an $E_8$-gauge field. 

The main purpose of this note is to prove integrality of the expression \eqref{witint} without the $E_8$-gauge field assumption, provided that we assume the \emph{stable} Hypothesis H of Fiorenza--Sati--Schreiber \cite{SS1}\cite[Section 3.8]{FSS}\cite[page 3]{BSS}. We also examine integrality of the expression \eqref{witint} assuming the \emph{unstable} Hypothesis H of Fiorenza--Sati--Schreiber \cite{Sati}\cite{FSS2}\cite{FSS3}\cite{HSS}\cite{SS}. The main theorems are given in section \ref{fluxq} (Theorem \ref{stablecoh} in the stable case and Theorem \ref{unstablecoh} in the unstable case).

\addtocontents{toc}{\SkipTocEntry}\subsection*{Acknowledgments} 
I would like to thank Urs Schreiber for suggesting this project as a follow-up to the results in \cite{GS} and for providing feedback on the first version of this article.

\section{Obstruction theory for stable $S^4$ revisited}

In this section, we revisit the construction of the first few stages of the Postnikov tower for $S^4$, as identified in \cite{GS}. We provide a self-contained proof of the identification of the $k$-invariants in the stable range. The overall strategy is as follows. We analyze the obstruction problem locally at each prime $p$ corresponding to a primary component of $\pi_{n+k}(S^{n+4})$ for $k\leq 11$ and then assemble the prime components to get the integral obstructions in this range. At each prime $p$, we move from one Postnikov section to the next by killing the first nonvanishing cohomology above degree $4$ at the previous stage. We do this by guessing what the $k$-invariant at the previous stage should be, taking the homotopy fiber of the proposed $k$-invariant, and showing that the relevant $\Z_p$-cohomology vanishes on the homotopy fiber. If this is the case, the homotopy fiber receives a map from $S^4$ which has higher connectivity (with respect to the $p$-primary component) than at the previous stage \footnote{This follows from the approximation theorem \cite[Theorem 4]{MT}}. We encourage the reader to refer to \cite[Chapter 12]{MT} for details on the general construction.

The homotopy groups of spheres have been computed in low degrees and we will use Toda's calculation in \cite{Toda} up to degree $11$ (See also \cite[page 339]{Hatcher} for a tabulation of Toda's calculation). These groups are displayed in the following table:

\begin{center}
\begin{tabular}{c|cccccccc}
$n$ & 4 & 5 & 6 & 7 & 8 & 9 & 10 & 11
\\
\hline
$\pi_n(S^4)$ & $\ZZ$ & $\Z_2$ & $\Z_2$ & $\Z\times \Z_{12}$ & $\Z_2^2$ & $\Z_2^2$ & $\Z_{24}\times \Z_3$ & $\Z_{15}$ 
\\
\hline
$\pi_n(\Sigma^{\infty}S^4)$ & $\Z$ & $\Z_2$ & $\Z_2$ & $\Z_{12}$ & 0 & 0 & $\Z_2$ & $\Z_{240}$

\end{tabular}.

\end{center}

From the table, we see that the primes corresponding to primary components of the stable homotopy groups $\pi_k(S^4)$ for $k\leq 11$ are $p=2, 3,5$. We analyze the obstruction theory at  each prime below. 
\begin{notation}\label{notation}
We use the following notation and conventions throughout
\begin{enumerate}
\item For an abelian group $A$ and a nonnegative integer $n$, we denote the $n$-fold shift of the corresponding Eilenberg--MacLane spectrum by $\Sigma^nHA$. 
\item For a spectrum $X$, we write $H^n(X;A):=[X,\Sigma^nHA]$, where the brackets denote the set of maps in the stable homotopy category.
\item As usual, we let ${\rm Sq}^n\in H^n(H\Z_2;\Z_2)$ denote the $n$-th Steenrod square. 
\item Let $p$ be prime. We let $\mathcal{P}_p^n\in H^{2n(p-1)}(H\Z_p;\Z_p)$ denote the $n$-th Steenrod power operation. 
\item We let $d_m\in H^{1}(H\Z_{p^m};\Z_p)$ denote the $m$-th Bockstein operation. 
\end{enumerate}
\end{notation}

\subsection{The prime $p=2$} We begin our analysis at the prime $p=2$. This analysis was carried out in \cite{MT} and we review the identification of the $k$-invariants here. At each stage of the tower, the first nonvanishing cohomology above degree 4 is finite and cyclic of order $2^m$, for some $m\geq 0$. At each stage, we choose a generator and kill the cohomology in this degree. This yields progressively better approximations for $S^4$. The discussion is summarized by the following proposition.

\begin{proposition}\label{Postnikov4}
The generators for the $\Z_2$-cohomology groups of the first five Postnikov sections of $\Sigma^{\infty}S^4\to \Sigma^4H\Z$, up to the first degree $n$ such that $H^n(X_i)\neq H^n(S^4)$, are given by the following table.
\begin{center}
\begin{tabular}{c|c|c|c|c|c|c}
\hline
 $n$ & $\Sigma^4H\Z$ & $X_1$ & $X_2$ & $X_3$ & $X_4$ & $X_5$
 \\
 \hline 
 0 & & & & & &
 \\
 1 & & & & & &
 \\
 2 & & & & & &
 \\
 3 & & & & & & 
 \\
 4 & $\rho_2\iota$ & $p^*\rho_2\iota$ & $p^*\rho_2\iota$ & $p^*\rho_2\iota$ & $p^*\rho_2\iota$ & $p^*\rho_2\iota$
 \\
 5 && & & &
 \\
 6 &${\rm Sq}^2\rho_2$ & & & &
 \\
 7 &    $\ast$          &  $\alpha_7$ & & &
 \\
 8 &  $\ast$          &    $\ast$      & $\mathfrak{Sq}^4$ & &
 \\
 9 &     $\ast$          &   $\ast$    &         $\ast$    & &
 \\
 10 & $\ast$ & $\ast$ & $\ast$& &
 \\
 11 &       $\ast$          &     $\ast$  &         $\ast$     & $p_{11}$ & 
 \\
 12 &       $\ast$          &    $\ast$   &        $\ast$     &   $\ast$ & $\mathfrak{Sq}^8$ 
\end{tabular} 
\end{center}
Each displayed entry generates a cyclic group of order 2, with the exception of $\mathfrak{Sq}^4$ and $\mathfrak{Sq}^8$. The class $\mathfrak{Sq}^4$ generates a cyclic group of order $8$ and $\mathfrak{Sq}^8$ generates a cyclic group of order $16$. Each generator is the $k$-invariant of the Postnikov section at the given stage.

\end{proposition}
\proof
Combine Propositions \ref{alpha7}, \ref{mod8thing}, \ref{p11}, \ref{mod16thing} and \ref{top}.
\endofproof 

Since $H^{5}(\Sigma^4H\ZZ; \Z_2)=0$, there is no obstruction in degree 5. The first nonvanishing obstruction occurs in degree $6$. We have 
$H^{6}(\Sigma^4H\ZZ; \Z_2)\cong \Z_2$, generated by ${\rm Sq}^2\rho_2$, where ${\rm Sq}^2$ is the 2nd Steenrod square and $\rho_2$ is the mod 2 reduction.
The first stage  of the Postnikov tower is the homotopy fiber sequence
$$
\xymatrix{
\Sigma^5H\ZZ_2 \ar[r]& X_1\ar[rr] \ar[d]^-{p} && \ast \ar[d]
\\
& \Sigma^4H\ZZ \ar[rr]^{{\rm Sq}^2\rho_2} && \Sigma^6H\ZZ_2.
}
$$
Since we are working stably, we can use the Serre cohomology sequence, i.e., the long exact sequence induced from the fiber/cofiber sequence $\Sigma^5H\ZZ_2\to X_1\to \Sigma^4H\ZZ$, to compute the cohomology of $X_1$. The sequence takes the form
\begin{equation}\label{serre}
\cdots \to H^*(\Sigma^4H\ZZ;\ZZ_2)\overset{p^*}{\to} H^*(X_1;\ZZ_2)\overset{i^*}{\to} H^*(\Sigma^5H\ZZ_2;\ZZ_2)\overset{\tau}{\to} H^{*+1}(\Sigma^4H\ZZ;\ZZ_2)\to \cdots 
\end{equation}
where $\tau$ is the transgression. To compute the cohomology of $X_1$, it is therefore sufficient to know the $\ZZ_2$ cohomology of $H\ZZ_2$, $H\ZZ$, and the transgression. The relevant cohomology is of course given in terms of the Steenrod algebra. For convenience, we recall the description of the cohomology of Eilenberg--MacLane spaces in terms of the Serre--Cartan basis in Appendix \ref{theapp}.

\begin{proposition}\label{alpha7}
We have $H^6(X_1;\Z_2)\cong 0$ and $H^7(X_1;\Z_2)\cong \Z_2$, generated by a class $\alpha_7$ such that $i^*(\alpha_7)={\rm Sq}^2$.
\end{proposition}
\proof
By construction, the fundamental class $\iota_5\in H^5(\Sigma^5H\Z_2;\Z_2)$ transgresses to ${\rm Sq}^2\rho_2 \in H^{6}(\Sigma^4H\ZZ; \Z_2)$. Since $H^6(\Sigma^5H\Z_2;\Z_2)\cong \Z_2$ is generated by ${\rm Sq}^1$ and the Steenrod operations commute with transgression, it follows from the Adem relation ${\rm Sq}^1{\rm Sq}^2={\rm Sq}^3$ that $\tau({\rm Sq}^1)={\rm Sq}^3\rho_2\neq 0$. From \eqref{serre}, we see that $H^6(X_1;\Z_2)\cong 0$, as claimed.

We also have $H^7(\Sigma^4H\ZZ;\ZZ_2)\cong 0$, $H^7(\Sigma^5H\Z_2;\Z_2)\cong \Z_2$, generated by ${\rm Sq}^2$, and $\tau({\rm Sq}^2)=0$. Therefore, $H^7(X_1;\Z_2)\cong \Z_2$, generated by a class $\alpha_7$ such that $i^*(\alpha_7)={\rm Sq}^2$. This completes the proof.
\endofproof

 We claim that $\alpha_7$ is the next obstruction. We consider the fiber sequence 
$$
\xymatrix{
\Sigma^6H\Z_2 \ar[r]& X_2\ar[rr] \ar[d]^-{p} && \ast \ar[d]
\\
& X_1 \ar[rr]^-{\alpha_7} && \Sigma^7H\Z_2\;.
}
$$

\begin{proposition}\label{mod8thing}
We have $H^7(X_2;\Z_2)\cong 0$ and $H^8(X_2;\Z_8)\cong \Z_8$, generated by a class $\mathfrak{Sq}^4$ whose mod 2 reduction is $\rho_2\mathfrak{Sq}^4=p^*{\rm Sq}^4\rho_2$, where $p:X_2\to \Sigma^4H\ZZ$ is the projection and $\rho_2$ is the mod 2 reduction.
\end{proposition}

\proof
We must compute the cohomology of $X_2$. This can again be done using the exact sequence
\begin{equation}\label{serre2}
H^*(X_1;\ZZ_2)\to H^*(X_2;\ZZ_2)\overset{i^*}{\to} H^*(\Sigma^6H\Z_2;\Z_2)\overset{\tau}{\to} H^{*+1}(X_1;\Z_2).
\end{equation}
Here, the transgressions are more involved, since they involve the cohomology of $X_1$ and not just $H\Z_2$. However, enough information can be extracted from $H^*(X_1;\Z_2)$ by the previous exact sequence to compute these. These  transgressions are computed in \cite[Lemma 1, Ch. 12]{MT} and the detailed calculation of $H^*(X_2;\Z_2)$ is given in \cite[Page 118-119]{MT}. The end result is that $H^7(X_2;\Z_2)\cong 0$ and $H^8(X_2;\Z_2)\cong \Z_2$, generated by $p^*{\rm Sq}^4\rho_2$, where $p:X_2\to \Sigma^4H\Z$ is the projection.  Since $d_1p^*{\rm Sq}^4=0$, where $d_1$ is the Bockstein (Notation \ref{notation}), it follows that $p^*{\rm Sq}^4$ is in the image of the mod 2-reduction $H\Z_4\to H\Z_2$. Choose a generator $s\in H^8(X_2;\Z_4)\cong \Z_4$ lifting $p^*{\rm Sq}^4$. Then also $d_2$ vanishes on $s$, by part 2. of \cite[Lemma 1]{MT}. Therefore, $p^*{\rm Sq}^4$ is the mod 2 reduction of a class $\mathfrak{Sq}^4\in H^8(X_2;\Z_8)\cong \Z_8$. 
\endofproof

It may seem odd that we have identified the $\Z_8$-cohomology in Proposition \ref{mod8thing} instead of the $\Z_2$-cohomology. However, as pointed out in \cite{MT}, the next obstruction in the Postnikov tower cannot be taken to be $p^*{\rm Sq}^4\rho_2\in H^8(X_2;\Z_2)$. The problem is that if we take $Y\to X_2$ to be the homotopy fiber of $p^*{\rm Sq}^4\rho_2$, then the fiber is $\Sigma^7H\Z_2$ has a class ${\rm Sq}^1\in H^8(\Sigma^7H\Z_2;\Z_2)$ that transgresses to $d_1 p^*{\rm Sq}^4=0$, where $d_1={\rm Sq}^1:H\Z_2\to \Sigma H\Z_2$ is the Bockstein. Hence $i^*:H^8(Y;\Z_2)\to H^8(\Sigma^7H\Z_2;\Z_2)$ must be nonzero by exactness, and we have not killed $H^8$. For this reason, we choose a generator $\mathfrak{Sq}^4\in H^8(X_2;\Z_8)$, given by Proposition \ref{mod8thing} and consider the fiber sequence
$$
\xymatrix{
\Sigma^7H\Z_8 \ar[r]& X_3\ar[rr] \ar[d]^-{p} && \ast \ar[d]
\\
& X_2 \ar[rr]^-{\mathfrak{ Sq}^4} && \Sigma^8H\Z_8\;.
}
$$

\begin{proposition}\label{p11}
We have $H^{8}(X_3;\Z_2)\cong H^9(X_3;\Z_2)\cong H^{10}(X_3;\Z_2)\cong 0$ and $H^{11}(X_3;\Z_2)\cong \Z_2$ generated by a class $p_{11}$ whose restriction to the fiber is ${\rm Sq}^4$. 
\end{proposition}
\proof
According to Proposition \ref{steenrod}, the cohomology of the fiber is $H^8(\Sigma^7H\ZZ_8;\Z_2)\cong \Z_2$, generated by $d_3$. The transgression is given by $\tau(d_3)=d_3p^*{\rm Sq}^4\iota_4\neq 0$ and the Serre cohomology sequence again implies that $H^8(X_3;\Z_2)\cong 0$, as desired. The cohomology $H^k(X_3;\Z_2)$ is given up to degree $k=10$ by \cite[Lemma 1, Ch. 12]{MT}, which gives the remaining identifications. 
\endofproof

Motivated by the previous proposition, we consider the homotopy fiber sequence
$$
\xymatrix{
\Sigma^{10}H\Z_2 \ar[r]& X_4 \ar[d]^-{p}\ar[rr] && \ast \ar[d]
\\
& X_3 \ar[rr]^-{p_{11}}  && \Sigma^{11}H\Z_2\;.
}
$$

\begin{proposition}\label{mod16thing}
We have $H^{11}(X_4;\Z_2)\cong 0$ and $H^{12}(X_4;\Z_{16})\cong \Z_{16}$ generated by a class $\mathfrak{Sq}^8$ whose mod 2 reduction is $\rho_2\mathfrak{Sq}^8=p^*{\rm Sq}^8\rho_2$. 
\end{proposition}
\proof
This follows from the identifications of the cohomology of $X_4$ in \cite[Lemma 1, Ch. 12]{MT}.
\endofproof

As in the case of $\mathfrak{Sq}^4$, we must work with a mod 16 refinement of the mod 2 class $p^*{\rm Sq}^8$. The final obstruction we consider is given by the class $\mathfrak{Sq}^8$ in Proposition \ref{mod16thing}. We consider the fiber sequence

$$
\xymatrix{
\Sigma^{11}H\Z_{16} \ar[r]& X_5 \ar[d]\ar[rr] && \ast \ar[d]
\\
& X_4 \ar[rr]^-{\mathfrak{Sq}^8}  && \Sigma^{12}H\Z_{16}\;.
}
$$
We have the following 
\begin{proposition}\label{top}
We have $H^{12}(X_5;\Z_2)\cong 0$.
\end{proposition}
\proof
This follows from the identification of the cohomology of $X_5$ in \cite[Lemma 1, Ch. 12]{MT}
\endofproof

\subsection{The primes $p=3,5$}

The main contribution to the Postnikov tower of $S^4$ is from the prime $p=2$ as we
saw above. To a lesser extent, the primes $p=3$ and $p=5$ also contribute. We summarize the relevant $\Z_3$-cohomology at each stage in the following proposition.

\begin{proposition} \label{Postnikov35}
Generators for the $\Z_3$-cohomology groups of the first two Postnikov sections $Y_1,Y_2$ of the map $\Sigma^{\infty}S^4\to \Sigma^4H\Z$ at the prime $p=3$, and the generator for the $\Z_5$-cohomology of the first Postnikov section $Z_1$ at the prime $5$ (up to the first degree nonvanishing degree $\geq 4$) are given by the following table.
\begin{center}
\begin{tabular}{c|c|c|c|c}
\hline
 $n$ & $\Sigma^4H\Z$ & $Y_1$ & $Y_2$ & $Z_1$
 \\
 \hline 
 0 & & & &
 \\
 1 & & & &
 \\
 2 & & & &
 \\
 3 & & & &
 \\
 4 & $\rho_3\iota$ & $p^*\rho_3\iota$ & $p^*\rho_3\iota$  & $p^*\rho_3\iota$ 
 \\
 5 && &&
 \\
 6 & &&
 \\
 7 &          &  &
 \\
 8 &  $\mathcal{P}^1\rho_3$         &    &
 \\
 9 &     $\ast$          &      &
 \\
 10 & $\ast$ & & &
 \\
 11 &       $\ast$          &    &
 \\
 12 &       $\ast$          &    $\beta_{12}$ & $\mathcal{P}^1_5\rho_5$
\end{tabular}
\end{center}
Each displayed entry generates a cyclic group of order 3. Each generator is the $k$-invariant of the Postnikov section at the given stage.
\end{proposition}

\proof
Combine Propositions \ref{p1}, \ref{vanishingstuff}, \ref{cubezero}, \ref{y2vanish}, and \ref{mod5}.
\endofproof
 
From the table in Lemma \ref{gensp}, we have  $H^8(\Sigma^4H\Z;\Z_3)\cong \Z_3$, generated by $\mathcal{P}_3^1\rho_3$,
\\
 $H^{12}(\Sigma^4H\Z;\Z_3)\cong \Z_3$, generated by $\mathcal{P}_3^2\rho_3$, and  $H^{12}(\Sigma^4H\Z;\Z_5)\cong \Z_5$, generated by $\mathcal{P}_5^1\rho_5$. From the above identifications, it readily follows that the first nonvanishing class is given by $\mathcal{P}^1_3\rho_3\in H^8(\Sigma^4H\Z,\Z_3)$. We claim that the first stage in the mod 3 Postnikov tower is given by the fiber sequence 
$$
\xymatrix{
\Sigma^7H\ZZ_3 \ar[r]& Y_1\ar[rr] \ar[d]^-{p} && \ast \ar[d]
\\
& \Sigma^4H\ZZ \ar[rr]^{\mathcal{P}^1_3\rho_3} && \Sigma^8H\ZZ_3.
}
$$
To prove this, we must show that $H^8$ has been killed on $Y_1$. As in the case $p=2$, we use the cohomology sequence
$$
\cdots \to H^*(\Sigma^4H\ZZ;\ZZ_3)\to H^*(Y_1;\ZZ_3)\overset{i^*}{\to} H^*(\Sigma^7H\ZZ_3;\ZZ_3)\overset{\tau}{\to} H^{*+1}(\Sigma^4H\ZZ;\ZZ_3)\to \cdots 
$$
where $\tau$ is the transgression. 

\begin{proposition}\label{p1}
We have $H^8(Y_1;\Z_3)\cong 0$. 
\end{proposition}
\proof We have $H^8(\Sigma^7H\ZZ_3;\ZZ_3)\cong \Z_3$, generated by $d_1$. Since $d_1$ commutes with transgression, it follows that $H^9(\Sigma^4H\Z;\Z_3)\cong \Z_3$ generated by $d_1\mathcal{P}^1_3\rho_3$. By construction, the fundamental class $\iota \in H^7(\Sigma^7H\Z_3;\Z_3)$ transgresses to $\mathcal{P}_3^1\rho_3$. Therefore, $d_1$ transgresses to $d_1\mathcal{P}_3^1\rho_3$. By exactness of the above sequence, this forces $H^8(Y_1;\Z_3)\cong 0$, as claimed. 
\endofproof

\begin{proposition}\label{vanishingstuff}
We have $H^9(Y_1;\Z_3)\cong H^{10}(Y_1;\Z_3)\cong H^{11}(Y_1;\Z_3)\cong 0$
\end{proposition}
\proof
Using the table in Appendix \ref{theapp}, we have exact sequences $0\to \Z_3(d_1)\overset{\tau}{\to} \Z_3(d_1\mathcal{P}^1_3)\to H^9(Y_1;\Z_3)\to 0$ and the transgression $\tau:d_1\mapsto d_1\mathcal{P}^1_3\iota_7$. Hence $H^9(Y_1;\Z_3)\cong 0$. We have the sequence $0\to H^{10}(Y_1;\Z_3)\to 0 \Rightarrow H^{10}(Y_1;\Z_3)\cong 0$ and $0\to H^{11}(Y_1;\Z_3)\to 0\Rightarrow H^{11}(Y_1;\Z_3)\cong 0$. 
\endofproof

\begin{proposition}\label{cubezero}
We have $H^{12}(Y_1;\Z_3)\cong \Z_3$, generated by a class $\beta_{12}$ whose restriction to the fiber is $\mathcal{P}^1_3d_1$. Moreover, $p^*\mathcal{P}^2_3\rho_3=0$ on $Y_1$. 
\end{proposition}
\proof
Again using the table in the appendix, we have the sequence 
$$0\to \Z_3(\mathcal{P}^1_3)\overset{\tau}{\to} \Z_3(\mathcal{P}^2_3\rho_3)\overset{p^*}{\to} H^{12}(Y_1;\Z_3)\to \Z_3(d_1\mathcal{P}^1_3,\mathcal{P}^1_3d_1)\overset{\tau}{\to} \Z_3(d_1\mathcal{P}^2_3\rho_3) $$
Since the map $\tau$ on the left is injective and both the domain and codomain are isomorphic to $\ZZ_3$, it follows that $\tau$ is surjective. Hence $\tau:\mathcal{P}^1_3\mapsto d_1\mathcal{P}^2_3\rho_3$ is an isomorphism and $p^*$ is the zero map. Hence $H^{12}(Y_1;\Z_3)$ injects into $\Z_3(d_1\mathcal{P}^1_3,\mathcal{P}^1_3d_1)$ and $p^*\mathcal{P}^2_3\rho_3=0$ on $Y_1$. The transgression maps $\tau:d_1\mathcal{P}^1_3\mapsto d_1\mathcal{P}^1_3\rho_3$. Hence $H^{12}(Y_1;\Z_3)\cong \Z_3(\mathcal{P}^1_3d_1)$, where the isomorphism is given by restriction to the fiber. We denote the corresponding generator of $H^{12}(Y_1;\Z_3)$ by $\beta_{12}$. This proves the claim.
\endofproof

Motivated by the above, we form the fiber sequence
$$
\xymatrix{
\Sigma^{11}H\Z_3\ar[r] & Y_2\ar[d] \ar[rr] && \ast\ar[d]
\\
& Y_1\ar[rr]^-{\beta{12}} && \Sigma^{12}H\Z_3
}
$$
It remains to prove that $\beta_{12}$ is indeed killed on $Y_2$. We use the fiber sequence 
\begin{equation}\label{seqys}
\cdots \to H^*(Y_1;\ZZ_3)\to H^*(Y_2;\ZZ_3)\overset{i^*}{\to} H^*(\Sigma^{11}H\ZZ_3;\ZZ_3)\overset{\tau}{\to} H^{*+1}(Y_1;\ZZ_3)\to \cdots 
\end{equation}
\begin{proposition}\label{y2vanish}
We have $H^{12}(Y_2;\Z_3)\cong 0$. 
\end{proposition}
\proof
Since $H^{12}(\Sigma^{11}H\Z_3;\Z_3)\cong \Z_3$, generated by $d_1$, the claim will follow provided that $\tau$ maps $d_1$ to a nontrivial element in $H^{13}(Y_1;\Z_3)$. To compute $H^{13}(Y_1;\Z_3)$, we use the sequence \eqref{seqys} to obtain a sequence 
$$
\Z_3(\beta_{12})\overset{i^*}{\to} \Z_3(\mathcal{P}^1_3d_1\rho_3,d_1\mathcal{P}^1_3)\overset{\tau}{\to} \Z_3(d_1\mathcal{P}^2_3)\to H^{13}(Y_1;\Z_3)\to \Z_3(d_1\mathcal{P}^1_3d_1)\overset{\tau}{\to} 0
$$
Now we have already seen that $\beta_{12}\mapsto \mathcal{P}^1_3d_1$ under $i^*$. By exactness, it follows that $\tau:d_1\mathcal{P}^1_3\mapsto d_1\mathcal{P}^2_3$. Again by exactness, it follows that $H^{13}(Y_1;\Z_3)\cong \Z_3$, generated by a class mapping to $d_1\mathcal{P}^1_3d_1$ under $i^*$. Finally, consider the pasting of homotopy pullback squares
$$
\xymatrix{
\Sigma^{11}H\Z_3\ar[r] &W\ar[d]\ar[r] & Y_2\ar[r]\ar[d] & \ast\ar[d]
\\
&\Sigma^7H\Z_3\ar[r]\ar@/_2pc/[rr]_{\mathcal{P}^1_3d_1} & Y_1\ar[r]^-{\beta_{12}} & \Sigma^{12}H\Z_3 
}
$$
where we have used that $i^*(\beta_{12})=\mathcal{P}^1_3d_1$. By construction, the fundamental class $\iota $ transgresses to $\mathcal{P}^1_3d_1$ along the map $W\to \Sigma^7H\Z_3$. Therefore, the transgression along $W\to \Sigma^7H\Z_3$ maps $\tau:d_1\to d_1\mathcal{P}^1_3d_1$. By the above commutativity, $d_1$ transgresses to a class $H^{13}(Y_1;\Z_3)$ whose restriction to the fiber is $d_1\mathcal{P}^1_3d_1$. It follows that $\tau(d_1)$ is nontrivial, as claimed. 
\endofproof

It only remains to identify the obstructions for the prime $p=5$. Here the first nontrivial class is $\mathcal{P}^1_5\rho_5\in H^{12}(\Sigma^4H\Z;\Z_5)$. We consider the Postnikov section 
$$
\xymatrix{
\Sigma^{11}H\Z_5\ar[r]&Z_1\ar[d]\ar[rr] && \ast\ar[d]
\\
 & \Sigma^4H\Z\ar[rr]^-{\mathcal{P}^1_5\rho_5} && \Sigma^{12}H\Z_5 .
}
$$

\begin{proposition}\label{mod5}
We have $H^{12}(Z_1;\Z_5)\cong 0$
\end{proposition}
\proof
By construction the fundamental class $\iota$ transgresses to $\mathcal{P}^1_5\rho_5$. The result is immediate from the Serre exact sequence. 
\endofproof

\subsection{Assembly into integral obstructions}

We are now ready to assemble the analysis at the primes $p=2,3,5$ into the integral obstruction theory for $S^4$. We have the following proposition.

\begin{proposition}\label{poststable} The integral Postnikov tower for the stable $4$-sphere takes the form

$$
\xymatrix{
 & \Sigma^{\infty}S^4\ar[d]&
 \\
 &\vdots \ar[d] &
 \\
 &W_5 \ar[d]&
 \\
\Sigma^{10}H\Z_2\ar[r] & W_4 \ar[d]\ar[rr]^-{(\mathfrak{Sq}^8,\beta_{12},p^*\mathcal{P}^1_5)} &&\Sigma^{12}H\Z_{16}\times \Sigma^{12}H\Z_3\times \Sigma^{12}H\Z_5
\\
 \Sigma^{7}H\Z_{24} \ar[r] & W_3\ar[d]\ar[rr]^-{p_{11}} && \Sigma^{11}H\Z_2
\\
   \Sigma^6H\Z_2 \ar[r] & W_2 \ar[d]\ar[rr]^-{(\mathfrak{Sq}^4,p^*\mathcal{P}^1_3)}  &&   \Sigma^8H\Z_8\times \Sigma^8H\Z_3
\\
  \Sigma^5H\Z_2\ar[r] & W_1 \ar[d]\ar[rr]^-{\alpha_7}  &&  \Sigma^7H\Z_2
   \\
     & \Sigma^4H\Z \ar[rr]^-{{\rm Sq}^2} && \Sigma^6H\Z_2
}
$$
where $p:W_i\to \Sigma^4H\Z$ is the projection. The $k$-invariants generate the cohomology of the corresponding Postnikov stage and satisfy the following conditions:
\begin{enumerate}
\item $\alpha_7$ restricts to ${\rm Sq}^2$ and $p_{11}$ restricts to ${\rm Sq}^4$ on the fiber;
\item Under mod 2 reduction, we have $\rho_2\mathfrak{Sq}^4=p^*{\rm Sq}^4$ and $\rho_2\mathfrak{Sq}^8=p^*{\rm Sq}^8$;
\item $\beta_{12}$ restricts to $\mathcal{P}^1_3d_1$ on the fiber.
\end{enumerate}
Moreover, we have the additional condition
\begin{enumerate}
\item[(4)] $p^*\mathcal{P}^2_3\rho_3=0$ on $W_3$.
\end{enumerate}

\end{proposition}
\proof
The first two stages of the tower coincide with the tower for the prime $p=2$, since the canonical map $S^4\to X_2$ induces an isomorphism on $2$-primary components of homotopy groups $\pi_n$ for $4<n\leq 6$ by construction, and the homotopy groups of $S^4$ are $2$-groups in this range. Set $W_1=X_1$ and $W_2=X_2$. 
Let $p:X_2\to \Sigma^4H\Z$ and $p_3:X_3\to X_2$ denote the projections in the Postnikov tower at the prime $2$ and let $q:Y_1\to \Sigma^4H\Z$ denote the projection in the Postnikov tower at the prime 3. Using  the pasting lemma and the two fiber sequences for the prime $2$ and $3$, we form iterated homotopy pullback squares
\begin{equation}\label{therealpostcomb}
\xymatrix{
\Sigma^7H\Z_8\times \Sigma^7H\Z_3\ar[r]\ar[d] & W_3\ar[r]\ar[d] & X_3\times Y_1\ar[rr]\ar[d]^-{p_3\times q_1} && \ast\ar[d]
\\
\ast\ar[r] & X_2\ar[r]^-{({\rm id},p_2)} & X_2\times \Sigma^4H\ZZ\ar[rr]^-{\mathfrak{Sq}^4\times \mathcal{P}^1_3} && \Sigma^8H\Z_8\times \Sigma^8H\Z_3.
}
\end{equation}
By the long exact sequence for the fiber sequence given by the left square, we have 
\begin{enumerate}
\item $\pi_i(W_3)\cong \pi_i(X_2)$ for $i\leq 6$ 
\item $\pi_7(W_3)\cong \Z_8\times \Z_3\cong \pi_7(X_3\times Y_1)$.
\end{enumerate}
 By construction of the Postnikov towers at the primes $2$ and $3$, we have maps $\Sigma^{\infty}S^4\to X_2$ and $\Sigma^{\infty}S^4\to X_3$ which induce isomorphisms on $2$-primary component of $\pi_i$ for $i\leq 6$ and $i\leq 7$, respectively. Similarly, we have a map $\Sigma^{\infty}S^4\to Y_1$ which induces an isomorphism on the $3$-primary component of $\pi_i$ for $i\leq 7$. 
The universal property of the fiber product in the center square of \eqref{therealpostcomb} gives a map $\Sigma^{\infty}S^4\to W_3$. Since $\pi_i(\Sigma^{\infty}S^4)$ has only $2$ and $3$-primary components for $i\leq 7$, it follows from the isomorphisms (1) and (2) above that the map $\Sigma^{\infty}S^4\to W_3$ induces an isomorphism on $\pi_i$ for $i\leq 7$. 

For the next stage of the tower, we proceed analogously as follows. Consider the pasting diagram 
$$
\xymatrix{
\Sigma^{10}H\Z_2 \ar[r]\ar[d] & W_4\ar[r]\ar[d] & X_4\ar[rr]\ar[d]^-{p_4} && \ast\ar[d]
\\
\ast\ar[r] & W_3\ar[r]^-{\rho_2} & X_3 \ar[rr]^-{p_{11}} && \Sigma^{11}H\Z_2.
}
$$
where $\rho_2:W_3\to X_3$ is the projection onto the second factor of the top center map in \eqref{therealpostcomb}. The universal property again produces a map $\Sigma^{\infty}S^4\to W_4$. Since $\pi_i(\Sigma^{\infty}S^4)$ has only 2 and 3 primary components for $4<i\leq 10$, and the map $\Sigma^{\infty}S^4\to X_4$ induces an isomorphism on 2 primary components in this range, the universal map induces an isomorphism on $\pi_i$ for $i\leq 10$. Finally, $W_5$ fits into a pasting diagram 
$$
\xymatrix{
\Sigma^{10}H\Z_{240}\ar[r]\ar[d] & W_5\ar[r]\ar[d] & X_5\times Y_2\times Z_1\ar[rr]\ar[d]^-{p_5\times q_2\times r_1} && \ast\ar[d]
\\
\ast\ar[r] & W_4\ar[r]^-{(\rho_2,\rho_3,p)} & X_4\times Y_1 \times \Sigma^4H\Z \ar[rr]^-{\mathfrak{Sq}^8\times \beta_{12}\times  \mathcal{P}^1_5} && \Sigma^{10}H\Z_{240}.
}
$$
where $\rho_2:W_4\to X_4$, and $\rho_3$ are again given by taking the top center maps above and projecting onto the corresponding factor. The same argument as before produces a map $\Sigma^{\infty}S^4\to W_5$ which induces an isomorphism on $\pi_i$ for $i\leq 11$. 

%
%
%
\endofproof

\section{The unstable obstructions}

So far, we have worked out the obstructions for lifting to the stable $4$-sphere. We now turn our attention to the unstable case. Here, the situation is more complicated, since we can no longer rely on the Serre cohomology sequence. Instead, we must use the Serre spectral sequence to compute the cohomology. We again examine the tower at each prime $p$ and assemble to an integral tower. 

The unstable homotopy groups of $S^4$ have been computed in low degrees and the relevant primes are $p=2,3,5$. 

\subsection{The prime $p=2$}

\begin{figure}
$$
\xymatrix@R=.2cm@C=.2cm{
q &\ar@{-}[dddddddddddd] &&&&&&&&&&&&&&&
\\
11 && \bullet\bullet\bullet\ar[dddddddddddrrrrrrrrrrrr]  &&&&  &&  &&&&&&&
\\
10 && {\color{red}{\bullet}}\bullet \ar[dddrrrr] \ar[ddddddddddrrrrrrrrrrr]  &&&& &&  &&&&&&&
\\
9 && {\color{red}{\bullet}} \bullet \ar[dddddddddrrrrrrrrrr]  &&&&  &&  &&&&&&&
\\
8 && \bullet\bullet \ar[ddddddddrrrrrrrrr] &&&&  &&  &&&&&&&
\\
7 && \bullet \ar[dddddddrrrrrrrr] &&&& \bullet\ar[drr]\ar[ddrrr]\ar[dddddddrrrrrrrr] &&  &&&&&&& 
\\
6 && \color{red}{\bullet}\ar[ddddddrrrrrrr] &&&& \color{red}{\bullet}\ar[ddddddrrrrrrr]\ar[drr] && \bullet & &&&&&&
\\
5 &&  \color{red}{\bullet} \ar[dddddrrrrrr] &&&& \color{red}{\bullet} \ar[dddddrrrrrr] && \bullet\ar[dddddrrrrrr] &\bullet&&&&& 
\\
4 &&&&&&&&&&&&
\\
3 && &&&&&&&&&&
\\
2 && &&&&&&&&&&
\\
1 &&&&&&&&&&&&
\\
0 && \bullet &&&& \bullet && \color{red}{\bullet}  & \color{red}{\bullet} & \bullet & & \color{red}{\bullet\bullet} & {\color{red}{\bullet}} \bullet & \bullet \bullet&
\\
\ar@{-}[rrrrrrrrrrrrrrr]&&&&&&&&&&&&&&&
\\
 && 0 & 1 & 2 & 3 & 4 & 5 & 6 & 7 & 8 & 9 & 10 & 11 & 12 &  p
}
$$
\caption{The differentials in the spectral sequence for the fibration $K(\Z_2,5)\to X_1\to K(\Z,4)$ up to the $E_{12}$-page. A dot corresponds to a generator in that degree. Black dots represent generators that are either in the kernel or cokernel of a differential and contribute to the cohomology.}\label{spectralseq}
\end{figure}
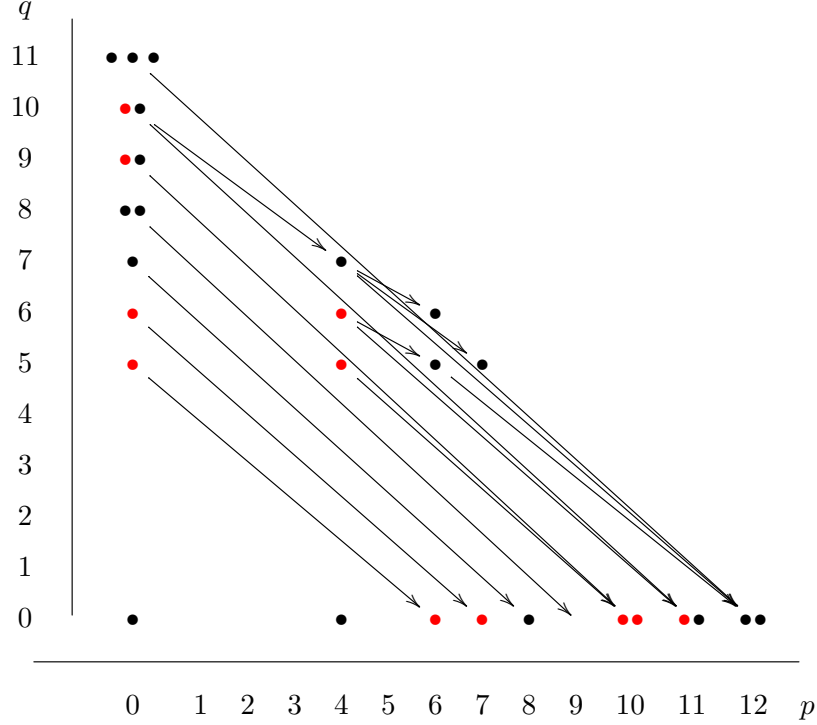

We begin with the most complicated case, when $p=2$. Since $H^5(K(\Z,4);\Z_2)=0$, there is no obstruction in degree $5$. The first nonvanishing obstruction occurs in degree $6$. We have $H^6(K(\Z,4);\Z_2)\cong \Z_2$, generated by ${\rm Sq}^2\iota_4$. The first stage is the same as in the stable case and is given by the fiber sequence 
$$
\xymatrix{
K(\Z_2,5) \ar[r]& X_1\ar[rr] \ar[d] && \ast \ar[d]
\\
& K(\Z,4)\ar[rr]^{{\rm Sq}^2\rho_2} && K(\Z_2,6)
}
$$
kills the class ${\rm Sq}^2\rho_2\iota_4$. Next, we must check that $H^6(X_1;\Z_2)=0$ and we will need the cohomology $H^*(X_1;\Z_2)$ in degrees $\leq 12$. We have the following
\begin{proposition}
The space $X_1$ has the following cohomology in degrees $n\leq 12$, where each displayed entry generates a copy of $\Z_2$.
\begin{center}
\hspace{-.2in}
\begin{tabular}{c|c|c|c}
\hline
$n$ & $H^{4+n}(K(\Z,4))$ generators & $H^{4+n}(K(\Z_2,5))$ generators & $H^{4+n}(X_1)$
\\
\hline 
0 & $\iota_4$ & & $p^*\iota_4$
\\
1 & & $\iota_5$ &
\\
2 & ${\rm Sq}^2\rho_2\iota_4$ & ${\rm Sq}^1\iota_5$  &
\\
3 & ${\rm Sq}^3\rho_2\iota_4$ & ${\rm Sq}^2\iota_5$ & $\alpha_7$
\\
4 & $\rho_2\iota_4^2$ & ${\rm Sq}^3\iota_5$, ${\rm Sq}^2{\rm Sq}^1\iota_5$ & $\beta_8$, $\gamma_8$, $p^*\rho_2\iota^2_4$
\\
5 & &${\rm Sq}^3{\rm Sq}^1\iota_5, {\rm Sq}^4\iota_5$ & $\delta_{9}$
\\
6 & ${\rm Sq}^4{\rm Sq}^2\rho_2\iota_4, {\rm Sq}^2\rho_2\iota_4\cdot \rho_2\iota_4$ & ${\rm Sq}^4{\rm Sq}^1\iota_5, \iota_5^2$ & $\epsilon_{10}$
\\
7 & ${\rm Sq}^5{\rm Sq}^2\rho_2\iota_4, {\rm Sq}^3\rho_2\iota_4\cdot \rho_2\iota_4$  &  ${\rm Sq}^5{\rm Sq}^1\iota_5, {\rm Sq}^4{\rm Sq}^2\iota_5, {\rm Sq}^1\iota_5\cdot \iota_5$
\\
8 & $\rho_2\iota_4^3, ({\rm Sq}^2\rho_2\iota_4)^2$ & ${\rm Sq}^4{\rm Sq}^2{\rm Sq}^1\iota_5$, ${\rm Sq}^4{\rm Sq}^2{\rm Sq}^1\iota_5, {\rm Sq}^5{\rm Sq}^2\iota_5, {\rm Sq}^2\iota_5\cdot \iota_5$
\\
\end{tabular}

\end{center}

\end{proposition}
\proof
The cohomology groups of the first 2 columns are computed in Lemma \ref{k4steenrod}. We use the Serre spectral sequence for the fibration $K(\Z_2,5)\to X_1\to K(\Z,4)$. The spectral sequence with differentials up to the $E_{12}$-page is displayed in Figure \ref{spectralseq}. 

 In degree $6$, we have two classes that can contribute, ${\rm Sq}^2\rho_2\iota_4$ and ${\rm Sq}^1\iota_5$. The transgression $\tau:\iota_5\mapsto {\rm Sq}^2\rho_2\iota_4$ and by the Adem relation ${\rm Sq}^1{\rm Sq}^2={\rm Sq}^3$, $\tau:{\rm Sq}^1\iota_5\mapsto {\rm Sq}^3\rho_2\iota_4$. Hence, both are killed. In degree 7, ${\rm Sq}^2\iota_5$ can contribute. Since $\tau:\iota_5\mapsto {\rm Sq}^2\rho_2\iota_4$, we have $\tau:{\rm Sq}^2\iota_5\mapsto {\rm Sq}^2{\rm Sq}^2\rho_2\iota_4={\rm Sq}^3{\rm Sq}^1\rho_2\iota_4=0$. Hence, ${\rm Sq}^2\iota_5$ survives to the $E_{\infty}$-page and $H^7(X_1;\Z_2)\cong \Z_2$, generated by a class $\alpha_7$ whose restriction to the fiber is ${\rm Sq}^2\iota_5$, as in the stable case. 
 
 Next, we look at degree 8. In this case, all three classes ${\rm Sq}^3\iota_5$, ${\rm Sq}^2{\rm Sq}^1\iota_5$ and $\rho_2\iota_4^2$ survive. Hence, $H^8(X_1;\Z_2)\cong \Z_2^3$, generated by three classes $\beta_8$, $\gamma_8$, and $p^*\rho_2\iota_4^2$. The first two classes restrict to ${\rm Sq}^3\iota_5$ and ${\rm Sq}^2{\rm Sq}^1\iota_5$ on the fiber, respectively. 
 
 For degree 9, observe that since $d_6$ satisfies the graded Leibniz rule, we have $d_6(\iota_4\cdot \iota_5)=d_6(\iota_4)\cdot \iota_5  + \iota_4\cdot d_6(\iota_5)=\iota_4\cdot {\rm Sq}^2\rho_2\iota_4 $. Therefore, $\iota_4\cdot \iota_5$ is killed by the differential and does not contribute to the cohomology. On the other hand, the Adem relations imply that ${\rm Sq}^3{\rm Sq}^1\iota_5$ transgresses to ${\rm Sq}^3{\rm Sq}^3\rho_2\iota_4={\rm Sq}^5{\rm Sq}^1\rho_2\iota_4=0$. Hence ${\rm Sq}^3{\rm Sq}^1\iota_5$ contributes and we get a single generator $\delta_{9}$ of $H^9(X_1;\Z_2)\cong \Z_2$ that restricts to ${\rm Sq}^3{\rm Sq}^1\iota_5$. The class ${\rm Sq}^4\rho_2\iota_4$ transgresses to ${\rm Sq}^4{\rm Sq}^2\rho_2\iota_4$, so this class does not contribute. 
 
 The cohomology in degree 10 is a bit more complicated. First, we must identify the differential $d_2$ on $\rho_2\iota_4\cdot {\rm Sq}^1\iota_5$. By the Leibniz rule, we have $d_2(\rho_2\iota_4\cdot {\rm Sq}^1\iota_5)=\rho_2\iota_4\cdot d_2({\rm Sq}^1\iota_5)=0$. Hence $\rho_4\iota_4\cdot {\rm Sq}^1\iota_5$ survives to the $E_7$-page. Again using the Leibniz rule, we have $d_7(\rho_2\iota_4\cdot {\rm Sq}^1\iota_5)=\rho_2\iota_4\cdot d_7({\rm Sq}^1\iota_5)=\rho_2\iota_4\cdot {\rm Sq}^3\rho_2\iota_4\neq 0$. Hence $\rho_4\iota_4\cdot {\rm Sq}^1\iota_5$ is killed on the $E_7$-page. Now from the previous identifications, we see that there is only only two possible classes that can contribute to the cohomology in degree 10: ${\rm Sq}^4{\rm Sq}^1\iota_5$ and $\iota_5^2$. The transgressions maps $\tau:{\rm Sq}^4{\rm Sq}^1\iota_5\mapsto  {\rm Sq}^4{\rm Sq}^3\rho_2\iota_4$. By the Adem relations, ${\rm Sq}^4{\rm Sq}^3\rho_2\iota_4={\rm Sq}^5{\rm Sq}^2\rho_2\iota_4$. Hence ${\rm Sq}^4{\rm Sq}^1\iota_5$ is killed on the $E_{11}$-page. On the other hand, $\tau(\iota^2_5)=0$ by the Leibniz rule. It follows that $H^{10}(X_1;\Z_2)\cong \Z_2$ generated by a class $\epsilon_{10}$ restricting to $\iota_5^2$ on the fiber. 
 
 Now in degree 11, $d_2:\rho_2\iota_4\cdot {\rm Sq}^2\iota_5\mapsto 0$, $d_3:\rho_2\iota_4\cdot {\rm Sq}^2\iota_5\to 0$. Therefore, $\rho_2\iota_4\cdot {\rm Sq}^2\iota_5$ survives to the $E_{\infty}$-page. The black dot at $(6,5)$ corresponds to $\iota_5\cdot {\rm Sq}^3\rho_2\iota_4$ and this class also survives to the $E_{\infty}$-page. Finally, all three classes ${\rm Sq}^5{\rm Sq}^1\iota_5$,${\rm Sq}^4{\rm Sq}^2\iota_5$ and $\iota_5\cdot {\rm Sq}^1\iota_5$ also contribute. The graded group on the $E_{\infty}$-page is $E_{\infty}^{0,11}\oplus E^{7,4}_{\infty}\oplus E^{6,5}_{\infty}\oplus E^{11,0}\cong \Z_2^6$, generated by the 5 classes above.
 \endofproof
 
The previous proposition establishes in particular that $H^6(X_1;\Z_2)\cong 0$ and $H^7(X_1;\Z_2)\cong \Z_2$, generated by a class $\alpha_7$. For the next stage of the tower, we consider the fiber sequence
$$
\xymatrix{
K(\Z_2,6) \ar[r]& X_2\ar[rr] \ar[d]^-{p} && \ast \ar[d]
\\
& X_1 \ar[rr]^{\alpha_7} && K(\Z_2,7)
}
$$
To prove that this is a better approximation of $S^4$, we must show that the cohomology $H^7(X_2;\Z_2)$ vanishes. The following proposition identifies the cohomology of $X_2$ in low degrees.

\begin{figure}
$$\xymatrix@R=.2cm@C=.2cm{
q &\ar@{-}[dddddddddddd] &&&&&&&&&&&&&&&
\\
9 && \bullet \bullet\ar[dddddddddrrrrrrrrrr]  &&&& \bullet\bullet &&  &&&&&
\\
8 && {\color{red}\bullet} \ar[ddddddddrrrrrrrrr] &&&& \bullet &&  &&&&&
\\
7 && \color{red}{\bullet} \ar[dddddddrrrrrrrr] &&&& \bullet &&  &&&&& 
\\
6 && {\color{red}\bullet}\ar[ddddddrrrrrrr] &&&& \color{red}{\bullet}\ar[ddddddrrrrrrr] &&  &&&&&
\\
5 &&   &&&&   &&  &&&& & 
\\
4 &&&&&&&&&&&&
\\
3 && &&&&&&&&&&
\\
2 && &&&&&&&&&&
\\
1 &&&&&&&&&&&&
\\
0 && \bullet &&&& \bullet && & {\color{red}\bullet} & {\color{red}\bullet}\bullet \bullet & {\color{red}\bullet} & \bullet & \color{red}{\bullet}
\\
\ar@{-}[rrrrrrrrrrrrrr]&&&&&&&&&&&&&&&&
\\
 && 0 & 1 & 2 & 3 & 4 & 5 & 6 & 7 & 8 & 9 & 10 & 11 & p
}
$$
\caption{The Serre spectral sequence for the fibration $K(\Z_2,6)\to X_2\to X_1$ up to the $E_{10}$-page. A single dot represents a single generator and a double dot represents a pair of generators. Black dots represent generators that are either in the kernel or cokernel of a differential and contribute to the cohomology.}\label{spectralseq2}
\end{figure}
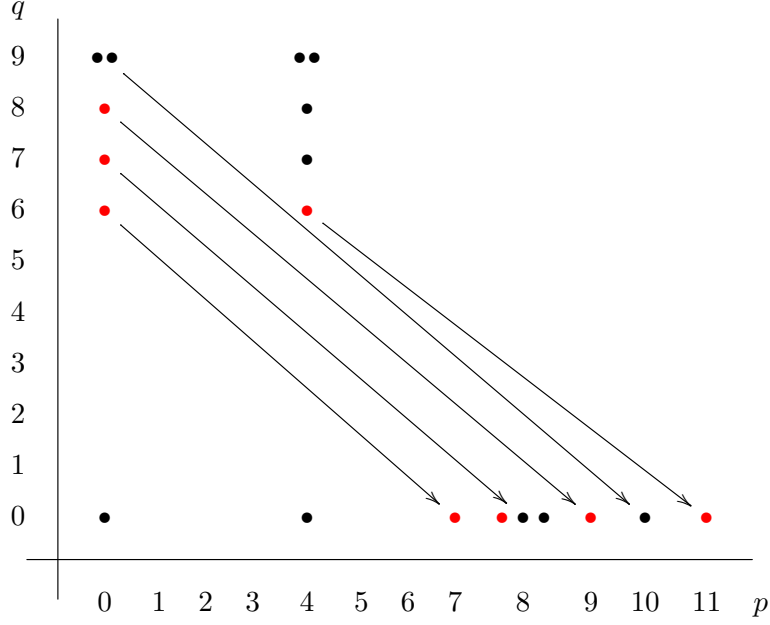

\begin{proposition}
The space $X_2$ has the following cohomology groups for $n\leq 12$, where each displayed entry generates a copy of $\Z_2$.
\begin{center}
\begin{tabular}{c|c|c|c}
\hline
$n$ & $H^{4+n}(X_1)$ generators & $H^{4+n}(K(\Z_2,6))$ generators & $H^{4+n}(X_2)$
\\
\hline 
0 & $p^*\iota_4$ & & $p^*\iota_4$
\\
1 & & &
\\
2 &  & $\iota_6$  &
\\
3 & $\alpha_7$ & ${\rm Sq}^1\iota_6$ & 
\\
4 & $\beta_8,\gamma_8, p^*\rho_4\iota_4^2$ & ${\rm Sq}^2\iota_6$  & $p^*\gamma_8$, $p^*\rho_2\iota^2_4$
\\
5 & $\delta_9$ &${\rm Sq}^2{\rm Sq}^1\iota_6, {\rm Sq}^3\iota_5$ & $\zeta_{9},\eta_9$
\\
6 & $\epsilon_{10}$ & ${\rm Sq}^3{\rm Sq}^1\iota_6, {\rm Sq}^4\iota_6$ & $p^*\epsilon_{10}$
\end{tabular}
\end{center}
\end{proposition}
\proof
The cohomology groups of the first column was computed in the previous proposition. The cohomology in the second column is computed in Lemma \ref{k4steenrod}. We use the Serre spectral sequence for the fibration $K(\Z_2,6)\to X_2\to X_1$.  The spectral sequence with differentials up to the $E_{10}$-page is displayed in Figure \ref{spectralseq2}.  

 In degree $7$, we have two classes that can contribute, ${\rm Sq}^1\iota_6$ and $\alpha_7$. By construction, the transgression $\tau:\iota_6\mapsto \alpha_7$ Therefore, $\tau:{\rm Sq}^1\iota_6\mapsto {\rm Sq}^1\alpha_7$. Since ${\rm Sq}^1\alpha_7$ restricts to ${\rm Sq}^1{\rm Sq}^2\iota_5={\rm Sq}^3\iota_5$ on the fiber, it follows that ${\rm Sq}^1\alpha_7\neq 0$. Write ${\rm Sq}^1\alpha_7=a\beta_8+b\gamma_8+cp^*\rho_2\iota_4$, $a,b,c\in \Z_2$. Let $j$ denote the fiber inclusion. Then 
 $${\rm Sq}^3\iota_6=j^*{\rm Sq}^1\alpha_7=a{\rm Sq}^3\iota_6+b{\rm Sq}^2{\rm Sq}^1\iota_6\Rightarrow a=1,b=0,$$
since ${\rm Sq}^3\iota_6$ and ${\rm Sq}^2{\rm Sq}^1\iota_6$ are independent. Hence, ${\rm Sq}^1\alpha_7=\beta_8$. This shows that $H^7$ is trivial. 

In $H^8$, we have two classes  four classes that can contribute, ${\rm Sq}^2\iota_6$, $\beta_8$, $\gamma_8$, and $p^*\rho_2\iota_4^2$. By the previous paragraph, $\tau:{\rm Sq}^1\iota_6\mapsto \beta_8$ is nontrivial. The cokernel is therefore isomorphic to $\Z^2_2$, generated by $\gamma_8$ and $p^*\rho_2\iota^2_4$. We have that ${\rm Sq}^2\iota_6$ transgresses to ${\rm Sq}^2\alpha_7$. This class restricts to ${\rm Sq}^2{\rm Sq}^2\iota_5={\rm Sq}^3{\rm Sq}^1\iota_5$ on the fiber. Since $\delta_9$ also restricts to ${\rm Sq}^2{\rm Sq}^2\iota_5={\rm Sq}^3{\rm Sq}^1\iota_5$, we must have $\delta_9={\rm Sq}^2\alpha_7$. Hence $\delta_9$ is killed. In total, we have $H^8(X_2;\Z_2)\cong \Z^2_2$, generated by $\gamma_8,p^*\rho_2\iota_4^2$. Now for $H^9$, the Leibniz rule implies $d_7:p^*\iota_4\cdot \iota_6\mapsto p^*\iota_4\cdot \alpha_7$. The two classes ${\rm Sq}^3\iota_6,{\rm Sq}^2{\rm Sq}^1\iota_6$ are the only classes that contribute. Both classes transgress to zero. Therefore, $H^9(X_2;\Z_2)\cong \Z_2^2$, generated by classes $\zeta_9$ and $\eta_9$ that restrict to ${\rm Sq}^3\iota_6,{\rm Sq}^2{\rm Sq}^1\iota_6$ (respectively), along the fiber inclusion.
 
\endofproof

Now to obtain the next obstruction, we must kill the cohomology in degree $8$. We have seen that there are two generators in degree 8. Leaving $p^*\rho_2\iota_4^2$ for a moment, observe that killing $p^*\gamma_8$ will not do the trick, since  the fiber will have a class ${\rm Sq}^1\iota_7$ transgressing to zero. However, observe that since the fiber inclusion $i$ is injective on $H^9$, we have
$$i^*d_1(p^*\gamma_8)=i^*p^*d_1(\gamma_8)=0\Rightarrow d_1(p_*\gamma_8)=0.$$
Therefore, $p^*\gamma_8$ is in the image of the mod 2 reduction $\Z_4\to \Z_2$. Let $\Gamma_8\in H^8(X_2;\Z_4)$ denote this class.  There is also the class $p^*\rho_2\iota_4$ and we must also kill this class. We form the fiber 
$$
\xymatrix{
K(\Z_4,7)\times K(\Z,7) \ar[r] & X_3\ar[d]\ar[rr] &&\ast\ar[d] 
\\
& X_2\ar[rr]^-{\Gamma_8,p^*\iota_4^2} && K(\Z_4,8)\times K(\Z,8)
}
$$ 
We claim that $H^8(X_3;\Z_2)\cong 0$.

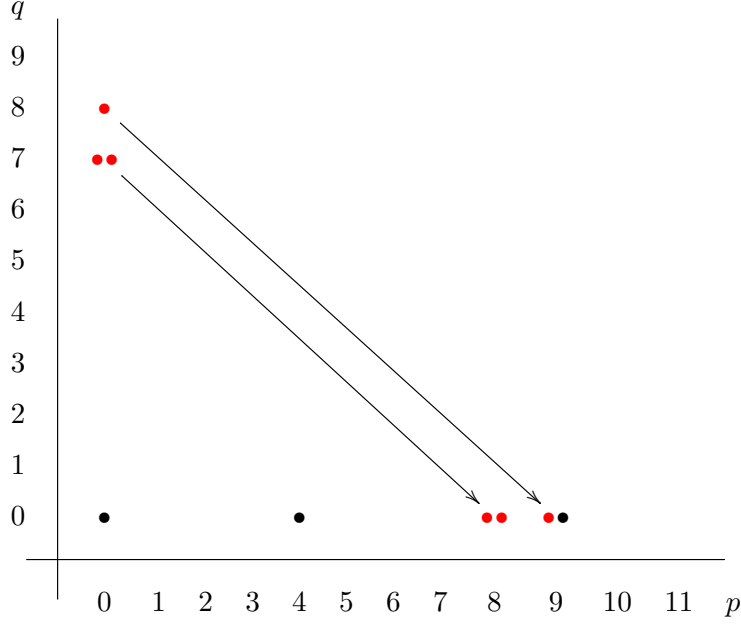
\begin{figure}
$$\xymatrix@R=.2cm@C=.2cm{
q &\ar@{-}[dddddddddddd] &&&&&&&&&&&&&&&
\\
9 &&   &&&& &&  &&&&&
\\
8 && {\color{red}\bullet  } \ar[ddddddddrrrrrrrrr] &&&&  &&  &&&&&
\\
7 && \color{red}{\bullet\bullet} \ar[dddddddrrrrrrrr] &&&& &&  &&&&& 
\\
6 && &&&& &&  &&&&&
\\
5 &&   &&&&   &&  &&&& & 
\\
4 &&&&&&&&&&&&
\\
3 && &&&&&&&&&&
\\
2 && &&&&&&&&&&
\\
1 &&&&&&&&&&&&
\\
0 && \bullet &&&& \bullet && && {\color{red}\bullet \bullet} & {\color{red}\bullet}\bullet & & 
\\
\ar@{-}[rrrrrrrrrrrrrr]&&&&&&&&&&&&&&&&
\\
 && 0 & 1 & 2 & 3 & 4 & 5 & 6 & 7 & 8 & 9 & 10 & 11 & p
}$$
\caption{The Serre spectral sequence for the fibration $K(\Z_4,7)\times K(\Z,7)\to X_3\to X_2$ up to the $E_0$-page. A single dot represents a single generator and a double dot represents a pair of generators. Black dots represent generators that are either in the kernel or cokernel of a differential and contribute to the cohomology.}\label{spectralseq3}
\end{figure}

\begin{proposition}
We have $H^8(X_3;\Z_2)\cong 0$.

\end{proposition}
\proof
Again using Lemma \ref{k4steenrod}. We see that the relevant portion of the Serre spectral sequence for the fibration $K(\Z_4,7)\times K(\Z,7)\to X_3\to X_2$ is given by Figure \ref{spectralseq3}, up to the $E_9$-page.
By construction $(\rho_4\iota_7,\iota_7)$ transgresses to $(\gamma_8,p^*\rho_2\iota_4^2)$. Now $d_1\gamma_8$ restricts along the fiber $K(\Z_2,5)\to X_1$ to ${\rm Sq}^1{\rm Sq}^2{\rm Sq}^1={\rm Sq}^3{\rm Sq}^1\iota_5\neq 0$. Since $\delta_9$ also restricts to this class, we deduce $d_1\gamma_8=\delta_9=\tau({\rm Sq}^2\iota_6)$. By the Bockstein lemma \cite[Theorem 1]{MT}, it follows that $d_2 \Gamma_8=d_1({\rm Sq}^2\iota_6)={\rm Sq}^1{\rm Sq}^2\iota_6={\rm Sq}^3\iota_6\neq 0$. Hence, $d_2\iota_7$ transgresses to $d_2\Gamma_8\neq 0$. Since $d_2\iota_7$ and $\Gamma_8$ are the only possible contributors to $H^8$, it follows that $H^8(X_3;\Z_2)\cong 0$, as claimed. 
\endofproof
\subsection{The prime $p=3$} 
At the prime 3, the first nonvanishing classes are $\mathcal{P}^1_3\rho_3\iota_4, \rho_3\iota_4^2\in H^8(K(\Z,4);\Z_3)$. We form the fiber sequence
$$
\xymatrix{
K(\Z_3,7)\times K(\Z,7)\ar[r] & Y_1\ar[rr]\ar[d] && \ast\ar[d]
\\
&K(\Z,4)\ar[rr]^-{\mathcal{P}^1_3\rho_3\iota_4, \rho_3\iota^2_4} && K(\Z_3,8)\times K(\Z,8)
}
$$

\begin{proposition}
We have $H^8(Y_1;\Z_3)\cong 0$
\end{proposition}
\proof
We will need the following cohomology groups, computed in Lemma \ref{k4steenrod}.
\begin{center}
\begin{tabular}{c|c|c}
\hline
$n$ & $H^{4+n}(K(\Z,4))$ generators & $H^{4+n}(K(\Z_3,7))$ generators 
\\
\hline 
0 & $\rho_3\iota_4$ & 
\\
1 & &  
\\
2 &  &  
\\
3 & & $\iota_7$ 
\\
4 & $\mathcal{P}^1_3\rho_3\iota_4, \rho_3\iota_4^2$ & $d_1\iota_7$ 
\\
5 & $d_1\mathcal{P}_3^1\rho_3\iota_4$ 
\\
6 &   & 
\\
7 & &  $\mathcal{P}^1_3\iota_7$
\\
8 & $\rho_3\iota_4^3, \mathcal{P}^1_3\rho_3\iota_4\cdot \rho_3\iota_4$ & $d_1\mathcal{P}^1_3\iota_7,\mathcal{P}^1_3d_1\iota_7$
\\
9 & $d_1\mathcal{P}^1_3\rho_3\iota_4\cdot \rho_3\iota_4$  
\\
10 &  & $\rho_3\iota_7^2$
\end{tabular}
\end{center}

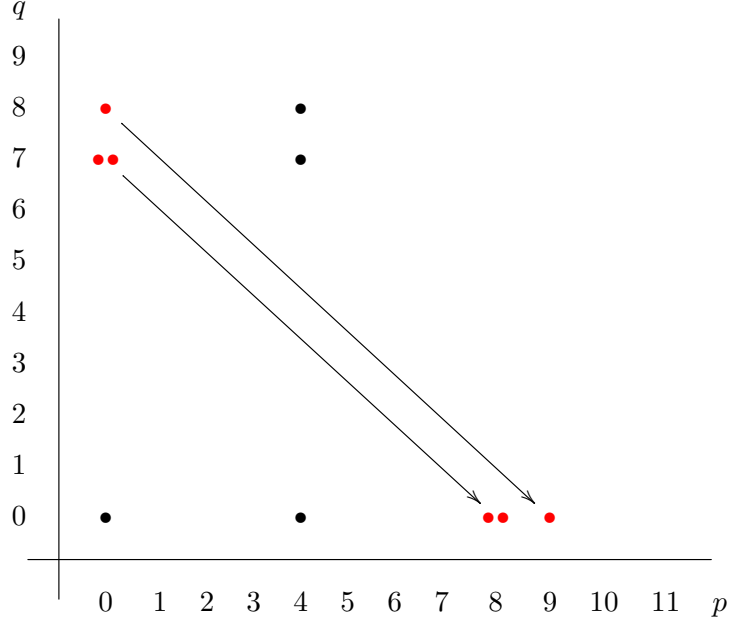
\begin{figure}
$$\xymatrix@R=.2cm@C=.2cm{
q &\ar@{-}[dddddddddddd] &&&&&&&&&&&&&&&
\\
9 && &&&&  &&  &&&&&
\\
8 && {\color{red}\bullet} \ar[ddddddddrrrrrrrrr] &&&& \bullet &&  &&&&&
\\
7 && \color{red}{\bullet\bullet } \ar[dddddddrrrrrrrr] &&&& \bullet &&  &&&&& 
\\
6 && &&&&  &&  &&&&&
\\
5 &&   &&&&   &&  &&&& & 
\\
4 &&&&&&&&&&&&
\\
3 && &&&&&&&&&&
\\
2 && &&&&&&&&&&
\\
1 &&&&&&&&&&&&
\\
0 && \bullet &&&& \bullet  && &  & {\color{red}\bullet \bullet }& {\color{red}\bullet} & & 
\\
\ar@{-}[rrrrrrrrrrrrrr]&&&&&&&&&&&&&&&&
\\
 && 0 & 1 & 2 & 3 & 4 & 5 & 6 & 7 & 8 & 9 & 10 & 11 & p
}$$
\caption{ The Serre spectral sequence for the fibration $K(\Z_3,7)\times K(\Z,7)\to Y_1\to K(\Z,4)$ up to the $E_8$-page. A single dot represents a single generator and a double dot represents a pair of generators. Black dots represent generators that are either in the kernel or cokernel of a differential and contribute to the cohomology.}\label{spectralseq4}
\end{figure}

 We use the Serre spectral sequence for the fibration $K(\Z_3,7)\times K(\Z,7)\to Y_1\to K(\Z,4)$. The spectral sequence is displayed in Figure \ref{spectralseq4}.
By construction $(\iota_7,\rho_3\iota_7)$ transgresses to $(\mathcal{P}^1_3\rho_3\iota_4,\rho_3\iota_4^2)$. Therefore, $d_1\iota_7$ transgresses to $d_1\mathcal{P}^1_3\rho_3\iota_4$. Hence $H^8(Y_1;\Z_3)$ is killed, as claimed.

\endofproof

\subsection{First three stages of the unstable tower.} We now assemble the previous obstructions to get the integral tower. We have the following proposition. 
pa
\begin{proposition}\label{postunstable} The first few layers of the Postnikov tower for the $4$-sphere are as follows
$$
\xymatrix{
 & S^4\ar[d]
 \\
 & \vdots \ar[d]&
 \\
 K(\Z_{11},8)\times K(\Z,7) \ar[r] & W_3\ar[d]^-{p_3} && 
\\
   K(\Z_{2},6) \ar[r] & W_2 \ar[d]^-{p_2}\ar[rr]^-{(\Gamma_8,p^*\mathcal{P}^1_3,\iota^2_4)}  &&   K(\Z_{4};8)\times K(\Z_3,8)\times K(\Z,8)
\\
  K(\Z_2,5)\ar[r] & W_1 \ar[d]^-{p_1}\ar[rr]^-{\alpha_7}  &&  K(\Z_2,7)
   \\
     & K(\Z,4) \ar[rr]^-{{\rm Sq}^2\rho_2\iota_4} && K(\Z_2,6)
}
$$
where the $k$-invariants are the classes identified in the previous section. The two unfamiliar $k$-invariants $\alpha_7$ and $\Gamma_8$ satisfy the following properties:

\begin{enumerate}
\item $\alpha_7$ restricts to ${\rm Sq}^2\iota_5$ on the fiber;
\item Under mod 2 reduction, we have that $\rho_2\Gamma_8=p_2^*\gamma_8$, where $\gamma_8\in H^8(W_1;\Z_2)$ restricts to ${\rm Sq}^2{\rm Sq}^1\iota_5$ on the fiber.
\end{enumerate}

\end{proposition}
\proof
The first two stages coincide with $X_1$ and $X_2$, which were constructed for the prime 2. We set $W_1=X_1$ and $W_2=X_2$. To get $W_3$, we proceed as follows. Let $p:X_2\to K(\Z,4)$ and $p_3:X_3\to X_2$ denote the projections in the Postnikov tower at the prime $2$ and let $q:Y_1\to K(\Z,4)$ denote the projection in the Postnikov tower at the prime 3. By successively using the pasting lemma and the two fiber sequences for the prime $2$ and $3$, we obtain a pasting diagram
$$
\xymatrix{
K(\Z_{12},7)\times K(\Z,7)\ar[r]\ar[d] & W_3\ar[r]\ar[d] & X_3\times Y_1\ar[rr]\ar[d]^-{p_3\times q_1} && \ast\ar[d]
\\
\ast\ar[r] & X_2\ar[r]^-{({\rm id},p_2)} & X_2\times K(\Z,4)\ar[rr]^-{\Gamma_8\times \mathcal{P}^1_3\rho_3\times \iota_4^2} && K(\Z_4,8)\times K(\Z_3,8)\times K(\Z,8).
}
$$
The same argument as in the stable case shows that the universal map $S^4\to W_3$ induces an isomorphism on $\pi_i$, for all $i\leq 8$. 
\endofproof

\section{Cohomology operations}

In this section, we describe the cohomology operations associated to the $k$-invariants of the Postnikov tower for the stable 4-sphere. In particular, we show that the operations obtained from the tower are additive modulo indeterminacy, whenever they are defined. This property will be used in the applications section to show that the integrality condition \eqref{witint} is the best possible, assuming the stable Hypothesis H.

There is a general prescription for producing cohomology operations from a Postnikov tower. Suppose we have a tower of fibrations along with $k$-invariants of the form 
\begin{equation}\label{towe.kinv}
\xymatrix{
 K_{n+k_2-1} \ar[r] & W_3\ar[d]&& 
\\
   K_{n+k_1-1} \ar[r] & W_2 \ar[d]\ar[rr]^-{\alpha_2}  &&   K_{n+k_2}
\\
 K_{n+k_0-1}\ar[r] & W_1 \ar[d]\ar[rr]^-{\alpha_1}  &&  K_{n+k_1}
   \\
     & K_n \ar[rr]^-{\alpha_0} && K_{n+k_0}
}
\end{equation}
where $K_n$ is an $n$-th Eilenberg-MacLane space. Set $A_i=\pi_i(K_i)$. The class of $\alpha_0$ defines a natural operation $\alpha_0:H^n(X;A_n)\to H^{n+k_0}(X;A_{n+k_0})$, by composing a map $f:X\to K_n$ with $\alpha_0$ and taking the homotopy class. 

The second obstruction $\alpha_1$ defines a secondary operation which is defined only on the kernel of $\alpha_0$. The operation is not uniquely defined and depends on a choice of lifting of a map $f:X\to K_n$ to $W_1$. The homotopy classes of lifts form a torsor for $H^{n+k_0-1}(X;A_{n+k_0})$. The operation obtained by restriction to the fiber $j^*\alpha_1:H^{n+k_1-1}(X;A_{n+k_0})\to H^{n+k_2}(X;A_{n+k_1})$ gives the indeterminacy of $\alpha_1$. In total, we have the secondary operation
$$
\alpha_1:\ker(\alpha_0)\to H^{n+k_1}(X;A_{n+k_1})/j^*\alpha_1H^{n+k_0-1}(X;A_{n+k_0})\;.
$$
The process continues with $\alpha_2$ being a tertiary obstruction defined on $\ker(\alpha_0)\cap \ker(\alpha_1)$ with indeterminacy given by the image of $i^*\alpha_2$.

\begin{proposition}\label{additive}

Modulo indeterminacy, the operations associated to the each of the obstructions
$\alpha_7,\mathfrak{Sq}^4, p_{11},\mathfrak{Sq}^8, \beta_{12}$
in Proposition \ref{poststable} are additive. That is, if $X$ is a spectrum and $x,y\in H^4(X;\Z)$ are such that the corresponding operation $\Phi$ is defined, then $\Phi$ is defined on $x+y$ and we have 
$$\Phi(x+y)=\Phi(x)+\Phi(y) \mod {\rm Ind}(\Phi).$$
\end{proposition}
\proof
This follows from stability of the operations. More generally, suppose we are given a tower of spectra of the form \eqref{towe.kinv} where each $\alpha_i:W_i\to K_{n+k_i}$ is a map into (a shift of) an Eilenberg--MacLane spectrum. Let $\Phi_i$ be the cohomology operation associated to $\alpha_i$. We claim that each $\Phi_i$ is additive modulo indeterminacy by induction on $i$. For $i=0$, the claim follows immediately, since $\alpha_0$ is a map of spectra (hence additive). Now assume the claim is true for $i$. Let $x,y:X\to K_n$ be such that $\Phi_{i+1}$ is defined on $x$ and $y$. Then there are lifts $\ell_x,\ell_y:X\to W_{i+1}$. Composing with the projection $W_{i+1}\to W_{i}$ yields lifts $\kappa_x,\kappa_y:X\to W_{i}$ such that $[\alpha_i(\kappa_x)]=\Phi_i(x)=0=\Phi_i(y)=[\alpha_i(\kappa_y)]$. By the induction hypothesis, $\Phi_i$ is defined for $x+y$. Hence, there is a lift $\kappa_{x+y}:X\to W_i$. Again by the induction hypothesis, $\Phi_i$ is additive modulo indeterminacy and it follows from the above that $[\alpha_i(\kappa_{x+y})]\in {\rm Im} j^*\alpha_i$, where $j:K_{n+k_{i-1}-1}\to W_i$ is the inclusion of the fiber. Since the set of lifts from $W_{i-1}$ to $W_i$ form a torsor for $H^{n+k_{i-1}-1}(X;A_{n+k_{i-1}})=[X,K_{n+k_{i-1}-1}]$, there is another lift $\kappa'_{x+y}:X\to W_i$ of $x+y$ such that $[\alpha_i(\kappa'_{x+y})]=0$. But then we can further lift $\kappa'_{x+y}$ through $W_{i+1}\to W_i$:
$$
\xymatrix{
& W_{i+1}\ar[d]
\\
X\ar[ru]^-{\ell_{x+y}}\ar[r]_-{\kappa'_{x+y}} & W_i\ar[r]^-{\alpha_i} & K_{n+k_i}.
}
$$
Hence, $\Phi_{i+1}$ is defined on $x+y$. Additivity is established as follows. We have a commutative diagram in the homotopy category of spectra
\begin{equation}\label{addcohopspec}
\xymatrix{
& W_{i+1}\vee W_{i+1}\ar[rr]^-{\alpha_{i+1}\vee \alpha_{i+1}}\ar[d]^-{+} && K_{n+k_2}\vee K_{n+k_2}\ar[d]^-{+}
\\
X\ar[r]_-{\ell_{x}+\ell_y}\ar[ru]^-{\ell_x\vee \ell_y} & W_{i+1}\ar[rr]^-{\alpha_{i+1}} && K_{n+k_2}.
}
\end{equation}
Since both $\ell_x+\ell_y$ and $\ell_{x+y}$ are lifts of $x+y$, the commutative diagram \eqref{addcohopspec} gives  
\begin{align*}
\Phi_{i+1}(x)+\Phi_{i+1}(y) &= [\alpha_{i+1}(\ell_x)]+[\alpha_{i+1}(\ell_y)]
\\
&=[\alpha_{i+1}(\ell_x+\ell_y)]
\\
& \equiv [\alpha_{i+1}(\ell_{x+y})]=\Phi_{i+1}(x+y)  \mod {\rm Ind}(\Phi_{i+1}) .
\end{align*}

\endofproof

\section{Flux quantization of the $C$-field in stable cohomotopy}\label{fluxq}

In this section, we use the analysis of the Postnikov tower for the stable and unstable $4$-sphere to obtain integrality conditions on the expression \eqref{witint}. Throughout this section, we let $\pi^4_s$ denote stable 4-cohomotopy, i.e., the cohomology theory represented by the spectrum $\Sigma^{\infty}S^4$. We let $\pi^4$ denote unstable 4-cohomotopy, i.e., the functor from the homotopy category of spaces to sets represented by $S^4$.

We begin with the stable case.

\begin{proposition}\label{cube6}
Let $M$ be a closed oriented 12-manifold. Suppose that $x\in H^4(M;\Z)$ is a degree 4 class lifting to $\pi^4_s(M)$. Then 
$$x^3\equiv 0 \mod 6.$$
\end{proposition}
\proof
If $x$ lifts to $\pi^4_s$, then in particular, there is a lift $\ell:M\to W_3$ of $x$ to the third stage of the Postnikov tower in Proposition \ref{poststable}. Then we must have 
$$
0=\ell^*\rho_2\mathfrak{Sq}^4=\ell^*p^*{\rm Sq}^4={\rm Sq}^4x=x^2\mod 2. 
$$
By Proposition \ref{poststable}, we also have that the class $p^*\mathcal{P}^2_3=0$ on $W_3$. Hence, 
$$0=\ell^*p^*\mathcal{P}^2_3=\mathcal{P}^2_3x=x^3\mod 3.$$
Together, the two conditions imply that 
$x^3\equiv 0 \mod 6$, as claimed.
\endofproof

The analogous condition on the cube in the unstable case is as follows.
\begin{proposition}\label{cube}
Let $M$ be a closed oriented 12-manifold. Suppose that $x\in H^4(M;\Z)$ is a degree 4 class lifting to $\pi^4(M)$. Then 
$$x^3= 0.$$
\end{proposition}
 \proof
 If $x$ lifts to $\pi^4$, then in particular, there is a lift $\ell:M\to W_3$ of $x$ to the third stage of the Postnikov tower in Proposition \ref{postunstable}. Then we must have 
$$
0=\ell^*p^*\iota^2=x^2. 
$$
and hence $x^3=0$ as claimed. 
 \endofproof
 
The next proposition shows that that the divisibility by 6 in Proposition \ref{cube6} is the best possible.
\begin{proposition}\label{quatproj}
Fix an orientation on $\HH P^1$. Consider the corresponding oriented 12-dimensional manifold $W=\HH P^1\times \HH P^1\times \HH P^1$. Let $x=u+v+w\in H^4(W)\cong \ZZ^3$ be the element given by the sum of positive generators of $H^4$ of each factor, with respect to the given orientation. Then $x$ lifts to $\pi_s^4(W)$ (stable cohomotopy in degree 4) and 
$$\langle x^3,[W]\rangle=6.$$
\end{proposition}
\proof
By the K\"unneth formula, we see that $H^*(W;\Z)$ is concentrated in degrees $n=0,4,8,12$. Fix an orientation on $\HH P^1$ and let $u,v$ and $w$ denote the dual of the fundamental class of each factor. Then $u,v,w$ generate $H^4\cong \Z^3$, $uv,wv,uw$ generate $H^8\cong \Z^3$, and $uvw$ generates  $H^{12}\cong \Z$. 

We analyze what the obstruction theory gives us for the class $u+v+w$. The only potential nonvanishing invariants are $\mathfrak{Sq}^4$, $\mathcal{P}^1_3$, $\mathcal{P}^1_5$ and $\beta_{12}$ and $\mathfrak{Sq}^8$. Since all the odd degree cohomology vanishes, and there is no cohomology in degree $6$, the indeterminacy of each corresponding cohomology operation vanishes. By Proposition \ref{additive}, all these operations are additive. Since $u,v$ and $w$ are pullbacks of the generator $H^4(\HH P^1;\Z)\cong \Z$ along the three projections $\HH P^1\times \HH P^1\times \HH P^1\to \HH P^1$, the corresponding cohomology operations all vanish on $u,v$ and $w$ for degree reasons. Therefore, $x=u+v+w$ lifts to a class $y\in \pi_s^4(W)$, as claimed.

It remains to prove that $\langle x^3,[W]\rangle=6$. We have $x^3=(u+v+w)^3=6uvw$. Since $u,v$ and $w$ are dual to the fundamental class of $\HH P^1$, $uvw$ is dual to the induced fundamental class of the threefold product. Therefore, 
$$\langle x^3,[W]\rangle=6\langle uvw,[W]\rangle=6,$$
as claimed.

\endofproof

We conclude with the main results. Theorem \ref{stablecoh} shows that integrality of the expression \eqref{witint} can be obtained assuming flux quantization of the 11D supergravity $C$-field in stable 4-cohomotopy, even without introducing a ``1-loop term" and an E8 gauge field. We also show that divisibility by $6$ is in general the best that one can do under this assumption. Theorem \ref{unstablecoh} shows that if we assume flux quantization in \emph{unstable} 4-cohomotopy, then the expression \eqref{witint} necessarily vanishes.

\begin{theorem}\label{stablecoh}
Let $M$ be a closed oriented 12-manifold. Let $v:\pi_s^4(M)\to H^4(M;\Z)$ denote the canonical map induced by a choice of generator $v\in H^4(S^4;\Z)\cong \Z$. Consider the function 
$$\pi_s^4(M)\to \ZZ, \quad x\mapsto \langle v(x)^3,[M]\rangle.$$
and let $(d_M)\subset \ZZ$ be the ideal generated by the image. Then $6\mid d_M$, moreover there exists a closed oriented 12-manifold $M$ such that $d_M=6$. 
\end{theorem}
\proof
The first statement follows Proposition \ref{cube6}. The second statement follows from  Proposition \ref{quatproj}, by taking $M$ to be three copies of $\HH P^1$.
\endofproof

\begin{theorem}\label{unstablecoh}
Let $M$ be a closed oriented 12-manifold. Let $v:\pi^4(M)\to H^4(M;\Z)$ denote the canonical map induced by a choice of generator $v\in H^4(S^4;\Z)\cong \Z$. Then the function 
$$\pi^4(M)\to \ZZ, \quad x\mapsto \langle v(x)^3,[M]\rangle$$
vanishes identically. 
\end{theorem}
\proof
This follows immediately from Proposition \ref{cube}.
\endofproof

\begin{appendices}
\section{The cohomology of Eilenberg--MacLane spaces}\label{theapp}

For $p$ prime, the structure of mod $p$ cohomology rings of Eilenberg-MacLane spaces was
determined by Cartan \cite{Cartan} and Serre \cite{Serre} in terms of admissible monomials of
 Steenrod reduced powers and the Bockstein. 
In this appendix, we recall these classical results.

\begin{definition}
Fix a prime $p$ and a multi-index $I:[r]\to \NN$. We call $I$ \emph{admissible for $p$} if the following hold:
\begin{itemize}
\item For all $j<r$, $i_j\geq pi_{j+1}$;
\item For all $j$, $i_j\equiv 0,1~{\rm mod}~2(p-1)$. 
\end{itemize}
\end{definition}

To describe the cohomology, we will need he \emph{excess} $e(I)$ of a sequence is defined as $e(I)=pi_1-(p-1)|I|$. Define $\epsilon(i_j)=i_j ~{\rm mod}~2(p-1)$ and consider also the multi-index $
\overline{I}=(\bar i_1,\hdots,\bar i_r)$
where
$$\bar i_j=\frac{i_j-\epsilon(i_j)}{2(p-1)}.$$

\begin{theorem}\label{steenrod}(Cartan)
Let $p=2$. Let $q$ be a nonnegative integer. Let $m$ be a positive integer. Then we have the following identifications.
\begin{itemize}
\item  $H^*(K(\ZZ,q),\Z_2)$ is a polynomial ring with generators ${\rm Sq}^I\rho_2\iota_q$, where $I$ runs through admissible sequences of excess $e(I)<q$ and where $i_r\neq 1$. 
\item The cohomology of $H^*(H(\ZZ_{2^m},q);\ZZ_2)$ is the polynomial ring with generators ${\rm Sq}^{I_m}$, where ${\rm Sq}^{I_m}:={\rm Sq}^I$ if $i_r\neq 1$ and ${\rm Sq}^{I_m}={\rm Sq}^{i_1}\cdots {\rm Sq}^{i_{r-1}}d_m$ if $i_r=1$, where $d_m$ is the $m$-th Bockstein, and $I$ runs through all admissible sequences of excess $e(I)<q$. 
\end{itemize}
\end{theorem}

\begin{theorem}\label{steenrodpower}(Cartan)
Let $p$ be prime and let $q$ be a nonnegative integer. Let $m$ be a positive integer. Then we have the following identifications.
\begin{itemize}
\item  $H^*(K(\ZZ,q),\Z_p)$ is a free graded commutative ring with generators in degrees $q+|I|$ given by $\mathcal{P}^{\bar I}\rho_3\iota_q$, if $i_r\neq 1$, and $d_1\mathcal{P}^{\bar I}\rho_3\iota_q$, if $i_r=1$, where $I$ runs through admissible sequences of excess $e(I)<(p-1)q$.  
\item The cohomology of $H^*(H(\ZZ_{p^m},q);\ZZ_p)$ is a free graded commutative ring with generators $\mathcal{P}^{I_m}$, where $\mathcal{P}^{I_m}:=\mathcal{P}^{\bar I}$ if $i_r\neq 1$ and $\mathcal{P}^{I_m}=\mathcal{P}^{\bar i_1}\cdots \mathcal{P}^{\bar i_{r-1}}d_m$ if $i_r=1$, where $d_m$ is the $m$-th Bockstein. 
\end{itemize}
\end{theorem}

\begin{lemma}\label{k4steenrod}
For $p=2$, we have the following admissible sequences $I$ and corresponding generators for $H^{|I|}(H\Z;\Z_2)$ and $H^{\vert I\vert}(H\Z_2;\Z_2)$, respectively.
\begin{center}
\begin{tabular}{c|c|c}
$H\Z$ 
\\
\hline
$|I|$ & $I$ & generators
\\
\hline 
0 & $(0)$ & $\iota$
\\
1 & $\emptyset$ &
\\
2 & $(2)$ & ${\rm Sq}^2\rho_2$
\\
3 & $(3)$ & ${\rm Sq}^3\rho_2$
\\
4 & $(4)$ & ${\rm Sq}^4\rho_2$
\\
5 & $(5)$ & ${\rm Sq}^5\rho_2$
\\
6 & $(4,2),(6)$ & ${\rm Sq}^4{\rm Sq}^2\rho_2, {\rm Sq}^6\rho_2$ 
\\
7 & $(5,2), (7)$ & ${\rm Sq}^5{\rm Sq}^2\rho_2, {\rm Sq}^7\rho_2$
\\
8 & $(6,2), (8)$ &  ${\rm Sq}^6{\rm Sq}^2\rho_2, {\rm Sq}^8\rho_2$
\\
9 & $(6,3), (7,2), (9)$ & ${\rm Sq}^6{\rm Sq}^3\rho_2, {\rm Sq}^7{\rm Sq}^2\rho_2, {\rm Sq}^9\rho_2$
\\
10 & $(7,3), (8,2), (10)$ & ${\rm Sq}^7{\rm Sq}^3\rho_2, {\rm Sq}^8{\rm Sq}^2\rho_2, {\rm Sq}^{10}\rho_2$
\end{tabular}
\end{center}

\begin{center}
\begin{tabular}{c|c|c}
$H\Z_2$ 
\\
\hline
$|I|$ & $I$ & generators
\\
\hline 
0 & $(0)$ & $\iota$
\\
1 & $(1)$ & ${\rm Sq}^1$
\\
2 & $(2)$ & ${\rm Sq}^2$
\\
3 & $(2,1), (3)$ & ${\rm Sq}^2{\rm Sq}^1, {\rm Sq}^3$
\\
4 & $(3,1), (4)$ & ${\rm Sq}^3{\rm Sq}^1, {\rm Sq}^4$
\\
5 & $(4,1), (5)$ & ${\rm Sq}^4{\rm Sq}^1, {\rm Sq}^5$
\\
6 & $(5,1), (4,2), (6)$ & ${\rm Sq}^5{\rm Sq}^1, {\rm Sq}^4{\rm Sq}^2, {\rm Sq}^6$
\\
7 & $(4,2,1), (6,1), (5,2), (7)$ & ${\rm Sq}^4{\rm Sq}^2{\rm Sq}^1, {\rm Sq}^6{\rm Sq}^1, {\rm Sq}^5{\rm Sq}^2, {\rm Sq}^7$
\end{tabular}
\end{center}
In the unstable case, taking into account the excess and products, we have the following admissible sequences and generators for the cohomology  $H^{4+|I|}(K(\Z,4);\Z_2)$, $H^{5+|I|}(K(\Z_2,5),\Z_2)$ and $H^{6+|I|}(K(\Z_2,6);\Z_2)$, respectively.
\begin{center}
\begin{tabular}{c|c|c}
$K(\Z,4)$ 
\\
\hline
$|I|$ & $I$ & generators
\\
\hline 
0 & $(0)$ & $\iota_4$
\\
1 & $\emptyset$ &
\\
2 & $(2)$ & ${\rm Sq}^2\rho_2\iota_4$
\\
3 & $(3)$ & ${\rm Sq}^3\rho_2\iota_4$
\\
4 & $\emptyset $ & $\rho_2\iota_4^2$
\\
5 & $\emptyset$ & 
\\
6 & $(4,2)$ & ${\rm Sq}^4{\rm Sq}^2\rho_2\iota_4, {\rm Sq}^2\rho_2\iota_4\cdot \rho_2\iota_4$ 
\\
7 & $(5,2)$ & ${\rm Sq}^5{\rm Sq}^2\rho_2\iota_4, {\rm Sq}^3\rho_2\iota_4\cdot \rho_2\iota_4$ 
\\
8 & $\emptyset$ & $\rho_2\iota_4^3, ({\rm Sq}^2\rho_2\iota_4)^2$
\\
9 & $(6,3)$ & ${\rm Sq}^6{\rm Sq}^3\rho_2, {\rm Sq}^3\rho_2\iota_4\cdot {\rm Sq}^2\rho_2\iota_4$
\\
10 & $\emptyset$ & ${\rm Sq}^2\rho_2\iota_4\cdot \rho_2\iota_4^2, ({\rm Sq}^3\rho_2\iota_4)^2$
\end{tabular}
\end{center}
\begin{center}
\begin{tabular}{c|c|c}
$K(\Z_2,5)$ 
\\
\hline
$|I|$ & $I$ & generators
\\
\hline 
0 & $(0)$ & $\iota_5$ 
\\
1 & $(1)$ & ${\rm Sq}^1\iota_5$  
\\
2 & $(2)$  & ${\rm Sq}^2\iota_5$ 
\\
3 & $(2,1), (3)$ & ${\rm Sq}^3\iota_5$, ${\rm Sq}^2{\rm Sq}^1\iota_5$ 
\\
4 & $(3,1), (4)$ &${\rm Sq}^3{\rm Sq}^1\iota_5, {\rm Sq}^4\iota_5$
\\
5 & $(4,1)$ &  ${\rm Sq}^4{\rm Sq}^1\iota_5, \iota_5^2$
\\
6 & $(5,1), (4,2)$ &   ${\rm Sq}^5{\rm Sq}^1\iota_5, {\rm Sq}^4{\rm Sq}^2\iota_5, {\rm Sq}^1\iota_5\cdot \iota_5$
\\
7 & $(4,2,1), (5,2)$ & ${\rm Sq}^4{\rm Sq}^2{\rm Sq}^1\iota_5, {\rm Sq}^5{\rm Sq}^2\iota_5, {\rm Sq}^2\iota_5\cdot \iota_5$
\end{tabular}
\end{center}
\begin{center}
\begin{tabular}{c|c|c}
$K(\Z_2,6)$ 
\\
\hline
$|I|$ & $I$ & generators
\\
\hline 
0 & $(0)$ & $\iota_6$ 
\\
1 & $(1)$ & ${\rm Sq}^1\iota_6$  
\\
2 & $(2)$  & ${\rm Sq}^2\iota_6$ 
\\
3 & $(2,1), (3)$ & ${\rm Sq}^3\iota_6$, ${\rm Sq}^2{\rm Sq}^1\iota_6$ 
\\
4 & $(3,1), (4)$ &${\rm Sq}^3{\rm Sq}^1\iota_6, {\rm Sq}^4\iota_6$
\\
5 & $(4,1), (5)$ &  ${\rm Sq}^4{\rm Sq}^1\iota_5, {\rm Sq}^5\iota_6$
\\
6 & $(5,1), (4,2)$ &   ${\rm Sq}^5{\rm Sq}^1\iota_6, {\rm Sq}^4{\rm Sq}^2\iota_6, \iota_6^2$
\\
7 & $(4,2,1), (6,1), (5,2)$ & ${\rm Sq}^4{\rm Sq}^2{\rm Sq}^1\iota_5, {\rm Sq}^5{\rm Sq}^2\iota_5, {\rm Sq}^6{\rm Sq}^1\iota_6, {\rm Sq}^1\iota_6\cdot \iota_6$
\end{tabular}
\end{center}
\end{lemma}
\proof
This is a direct application of Theorem \ref{steenrod}, where we take $q\to \infty$ for the first claim.
\endofproof

\begin{lemma}\label{gensp}
For $p=3$, we have the following admissible sequences $I$ and corresponding generators for $H^{|I|}(H\Z;\Z_3)$ and $H^{\vert I\vert}(H\Z_3;\Z_3)$, respectively.
\begin{center}
\begin{tabular}{c|c|c}
$H\Z$ 
\\
\hline
$|I|$ & $I$ & generators
\\
\hline 
0 & $(0)$ & $\iota$
\\
1 & $\emptyset$ &
\\
2 & $\emptyset$ &
\\
3 & $\emptyset$ &
\\
4 & $(4)$ & $\mathcal{P}^1_3\rho_3$
\\
5 & $(5)$ & $d_1\mathcal{P}^1_3\rho_3$
\\
6 & $\emptyset$ & 
\\
7 & $\emptyset$ & 
\\
8 & $(8)$ & $\mathcal{P}^2_3\rho_3$
\\
9 & $(9)$ & $d_1\mathcal{P}^2_3\rho_3$
\\
10 & $\emptyset$ &
\end{tabular}
\begin{tabular}{c|c|c}
$H\Z_3$ 
\\
\hline
$|I|$ & $I$ & generators
\\
\hline 
0 & $(0)$ & $\iota$
\\
1 & $(1)$ & $d_1$
\\
2 & $\emptyset$ &
\\
3 & $\emptyset$ &
\\
4 & $(4)$ & $\mathcal{P}^1_3$
\\
5 & $(5)$, $(4,1)$ & $d_1\mathcal{P}^1_3$, $\mathcal{P}^1_3d_1$
\\
6 & $(5,1)$ & $d_1\mathcal{P}^1_3d_1$
\\
7 & $\emptyset$
\\
 & 
\\
 & 
 \\
 &
\end{tabular}
\end{center}
In the unstable case, taking into account the excess and products, we have the following admissible sequences and generators for  $H^{4+|I|}(K(\Z,4);\Z_3)$ and $H^{7+|I|}(K(\Z_3,7),\Z_3)$, respectively.
\begin{center}
\begin{tabular}{c|c|c}
$K(\Z,4)$ 
\\
\hline
$|I|$ & $I$ & generators
\\
\hline 
0 & $(0)$ & $\rho_3\iota_4$
\\
1 & $\emptyset$ &
\\
2 & $\emptyset$ &
\\
3 & $\emptyset$ &
\\
4 & $(4)$ & $\mathcal{P}^1_3\rho_3\iota_4, \rho_3\iota_4^2$
\\
5 & $(5)$ & $d_1\mathcal{P}^1_3\rho_3\iota_4$
\\
6 & $\emptyset$ & 
\\
7 & $\emptyset$ & 
\\
8 & $\emptyset $ & $\rho_3\iota^3_4, \mathcal{P}^1_3\rho_3\iota_4\cdot \rho_3\iota_4$
\\
9 & $\emptyset $ & $d_1\mathcal{P}^1_3\rho_3\iota_4\cdot \rho_3\iota_4$
\\
10 & $\emptyset$ &
\end{tabular}
\begin{tabular}{c|c|c}
$K(\Z_3,7)$ 
\\
\hline
$|I|$ & $I$ & generator
\\
\hline 
0 & $(0)$ & $\iota_7$
\\
1 & $(1)$ & $d_1\iota_7$
\\
2 & $\emptyset$ &
\\
3 & $\emptyset$ &
\\
4 & $(4)$ & $\mathcal{P}^1_3\iota_7$
\\
5 & $(5)$, $(4,1)$ & $d_1\mathcal{P}^1_3\iota_7$, $\mathcal{P}^1_3d_1\iota_7$
\\
6 & $(5,1)$ & $d_1\mathcal{P}^1_3d_1\iota_7$
\\
7 & $\emptyset$ & $\rho_3\iota_7^2$
\\
 & 
\\
 & 
 \\
 &
\end{tabular}
\end{center}
\end{lemma}

  \end{appendices}

\end{document}